\documentclass[a4paper,11pt]{article}
\usepackage{graphicx}
\usepackage{latexsym,color}
\usepackage{amsmath,amsfonts,amssymb,amsthm}
\newcommand{\CB}{\Box}
\newtheorem{thm}{Theorem}[section]
\newtheorem{lem}[thm]{Lemma}

\newtheorem{cor}[thm]{Corollary}
\newtheorem{prop}[thm]{Proposition}
\newtheorem{conj}[thm]{Conjecture}
\newtheorem{rem}[thm]{Remark}

\input epsf

\title{The Stern Sequence and Moments of Minkowski's Question Mark Function}
\author{Roland Bacher}

\begin{document}
\maketitle


{\sl Abstract\footnote{
Math. class:  Primary: 11A55, Secondary: 11B57.}: 
We use properties of the Stern Sequence for numerical computations
of moments $\int_0^1 t^nd?(t)$ associated to Minkowski's Question 
Mark function.}
\vskip.5cm

\section{Introduction}

Minkowski's question mark function $x\longmapsto ?(x)$ and its 
inverse function, Conway's box function $x\longmapsto \CB(x)$, are 
related to continued fraction expansions, transcendence properties
and probabilistic distributions of rationals in the Calkin-Wilf tree.
Denjoy proved apparently that $?(x)$ is monotonic continuous and singular
(derivable on a set of full measure with zero derivative on this set),
see \cite{Denj}.
Using a functional equation satisfied by $?(x)$, Alkauskas investigated
the sequence $m_0,m_1,\dots,m_n=\int_0^1 x^nd?(s)$ 
of moments of the probability density $d?$ in a series of articles.
Denoting by $\CB(y)$ the reciprocal function, known as 
\emph{Conway's Box function}, 
of the increasing homeomorphism $?:[0,1]\longrightarrow [0,1]$,
the substitution $t=?(x)$ (with $d?(x)=dt$ and $x=\CB(t)$) yields
\begin{align}\label{defmoments}
m_n&=\int_0^1(\CB(t))^ndt\ .
\end{align}
In the present paper we link these moments to the Stern sequence
(which underlies the Calkin-Wilf tree)  $s(0)=0,s(1)=1,
s(2n)=s(n),s(2n+1)=s(n)+s(n+1),\ n\geq 1$. This gives new proofs for
many results of Alkauskas, see for example \cite{AlMC}, \cite{AlGMJ}, 
\cite{Alkasympt}. It also leads to the discovery of some
new properties. 

The sequel of the paper is organized as follows:

Section \ref{sectConw}
links the Stern sequence with Conway's Box function $\CB$ 
appearing in (\ref{defmoments}). 

Section \ref{sectmink?} recalls properties of Minkowski's question mark function.

Section \ref{sectusefulid} lists a few well-known identities among 
binomial coefficients and elements of the Stern sequence for later use.

Section \ref{sectsimpleset} 
presents a set of linear relations obtained by considering Riemann
sums for $\int_0^1 \CB(x)^ndx$.
These relations differ from the relations found by Alkauskas:
they are perhaps slightly simpler but more interestingly,
a crude spectral analysis of the underlying linear operator $T$ is
easy. $T$ has a unique eigenvector $(m_0,m_1,m_2,\dots)$ of eigenvalue $1$.
All other eigenvalues belong to the closed complex
disc of radius $1/2$. The maximal error of the associated algorithm
is thus roughly halved at each iteration.

Section \ref{sectformAlk} 
discusses a different set of Riemann sums which leads to 
linear relations used by Alkauskas.

We extend in Section \ref{sectholom} 
the moment-function $n\longmapsto 
m_n$ to an entire function $z\longmapsto m_z$ for $z\in \mathbb C$.

A computation of the derivative of this function at $0$ to high 
accuracy suggests the conjectural identities
$$\log 2=\sum_{n=1}^\infty \frac{m_n}{2n}\left(1+\frac{1}{2^{n-1}}\right)
=\sum_{n=1}^\infty \frac{m_n}{n}\left(\frac{1}{2^n}-(-1)^n\right)$$
given in Section \ref{sectconjrel}.

Section \ref{sectthirdset}
introduces a third type of Riemann sums, particularly well
suited for asymptotic computations. The resulting asymptotic formula
\begin{align}\label{asymptapprintro}
m_n&\sim \sum_{j=0}^\infty
\frac{(\log 2)^j}{j!}m_j\sum_{h=2}^\infty\frac{1}{2^h}\left(1-\frac{1}{h}
\right)^n
\end{align}
is the object of Section \ref{sectasympt}. 
It is more complicated but
experimentally more accurate than Alkauskas's asymptotic formula 
given in \cite{Alkasympt}. Alkauskas's formula can however be deduced from 
(\ref{asymptapprintro}) by a simple application of Laplace's method.

Section \ref{sectsecasympt} derives a second asymptotic formula related to 
(\ref{asymptapprintro}) by a finer subdivision in the underlying
Riemann sum. Since this should lead to slightly more accurate
results, we consider (admittedly in a not completely rigorous way) in
Section \ref{secterror} the difference between the two formulae as 
a measure of accuracy for (\ref{asymptapprintro}).

Section \ref{sectvalneg} is devoted to values $m_{-n}$ of moments 
at negative integers. This leads to a sequence of identities 
among $m_0,m_1,m_2,\dots$. The two initial identities are
$$\sum_{j=0}^\infty m_j=\frac{5}{2}\qquad\hbox{  and  }\qquad
\sum_{j=1}^\infty j\ m_j=m_2+\frac{11}{2}\ .$$

Finally, Section \ref{sectgeommeans} discusses the starting point of this
work: asymptotics for $\prod_{j=2^n}^{2^{n+1}}s_j$ allowing to compute
some geometric means for values of the Stern sequence.

\section{Conway's box function}\label{sectConw}

We denote by $\mathcal D=\mathbb Z[1/2]\cap [0,1]$ the subset of all 
rational dyadic numbers in $[0,1]$. The restriction 
to $\mathcal D$ of \emph{Conway's Box function}
$\CB$ is recursively defined as follows: $\CB(0)=0,
\CB(1)=1$ and
$$\CB\left(\frac{2m+1}{2^{n+1}}\right)=\frac{a+c}{b+d}$$
if $\CB\left(\frac{m}{2^n}\right)=\frac{a}{b}$
and $\CB\left(\frac{m+1}{2^n}\right)=\frac{c}{d}$
where $a,b$, respectively $c,d$, are coprime natural numbers.
The values $\CB(m/16)$ for $m=0,\dots,16$ are:
$$\begin{array}{c|ccccccccccccccccccccccccccc}
m&0&1&2&3&4&5&6&7&8&9&10&11&12&13&14&15&16\\
\hline
\CB\left(\frac{m}{16}\right)
&\frac{0}{1}&\frac{1}{5}&\frac{1}{4}&\frac{2}{7}&\frac{1}{3}&
\frac{3}{8}&\frac{2}{5}&\frac{3}{7}&\frac{1}{2}&\frac{4}{7}&\frac{3}{5}&
\frac{5}{8}&\frac{2}{3}&\frac{5}{7}&\frac{3}{4}&\frac{4}{5}&\frac{1}{1}
\end{array}$$

Values of $\CB$ for arguments in $\mathcal D$ are easy to compute as follows:
We define the {\it Stern-sequence} $s(0),s(1),s(2),\dots$
recursively by $s(0)=0,s(1)=1$,
$s(2n)=s(n)$ and $s(2n+1)=s(n)+s(n+1)$. Its first coefficients 
are given by 
$$\begin{array}{|c||c|c|c|c|c|c|c|c|c|c|c|c|c|c|c|c|}
\hline
n&0&1&2&3&4&5&6&7&8&9&10&11&12&13&14&15\\
\hline
s(n)&0&1&1&2&1&3&2&3&1&4&3&5&2&5&3&4\\
\hline\hline
n&16&17&18&19&20&21&22&23&24&25&26&27&28&29&30&31\\
\hline
s(n)&1&5&4&7&3&8&5&7&2&7&5&8&3&7&4&5\\
\hline
\end{array}$$
The main tool used in this paper is the following simple observation
which defines $\CB$ on $\mathcal D$ in terms of the Stern-sequence:
\begin{prop}\label{propfond}
We have $$\CB\left(\frac{m}{2^n}\right)=\frac{s(m)}{s(2^n+m)}$$
for all natural integers $m,n$ such that $0\leq m\leq 2^n$.
\end{prop}
We leave the easy proof to the reader.\hfill$\Box$

Since $\frac{a}{b}<\frac{c}{d}$ with $bd>0$ implies $\frac{a}{b}<\frac{a+c}{b+d}<
\frac{c}{d}$, the function $\CB$ is strictly increasing. 
Induction on $k$ shows
\begin{align}\label{formCkab}
\CB\left(\frac{m}{2^n}\pm \frac{1}{2^{n+k}}\right)&=\frac{ka+c}{kb+d},
\ k\geq 1
\end{align}
if $\CB\left(\frac{m}{2^n}\right)=\frac{a}{b}$ 
and $\CB\left(\frac{m}{2^n}\pm \frac{1}{2^n}\right)=\frac{c}{d}$
(with $a,b$ and $c,d$ pairs of coprime natural numbers). In
particular, $\CB$ extends to a strictly increasing continuous function 
(still denoted)
$\CB:[0,1]\longrightarrow[0,1]$. 
Since 
$$\lim_{k\rightarrow\infty}\frac{\CB
\left(\frac{m}{2^n}\pm \frac{1}{2^{n+k}}\right)- \CB
\left(\frac{m}{2^n}\right)}{\pm 2^{-n-k}}=\frac{2^{n+k}}{b(kb+d)}=\infty,$$
the function $\CB$ has a vertical tangent at dyadic arguments.


\begin{prop}\label{propsymC} We have 
\begin{align}\label{Csymmetric}
\CB(x)&=1-\CB(1-x)\ .
\end{align}
\end{prop}

\proof{} 
Continuity of Conway's 
Box function implies that 
it is enough to prove Proposition \ref{propsymC} for all 
dyadic rationals of the form $\frac{m}{2^n}$.
This is done by induction using the trivial identity 
$\frac{a+c}{b+d}=1-\frac{b-a+d-c}{b+d}$.
\endproof

\begin{cor}\label{corsymmetric}
The function $x\longmapsto 2\CB\left(x-\frac{1}{2}\right)-1$
is symmetric.
\end{cor}

Thus we have
\begin{align}\label{idparity}
0&=\int_0^1(\CB(x)-1/2)^{2n+1}dx=\sum_{k=0}^{2n+1}{2n+1\choose k}
(-2)^{k-2n-1}m_k
\end{align}
for every odd natural number $2n+1$. This can be restated as:

\begin{cor}\label{coridsymmn} For all $n\geq 0$ we have the identity
\begin{align}\label{idm2np1}
m_{2n+1}&=&\frac{1}{2^{2n+1}}\sum_{k=0}^{2n}(-2)^k{2n+1\choose k}m_k\ .
\end{align}
In particular, $m_{2n+1}$ is 
a $\mathbb Z\left[\frac{1}{2}\right]-$linear 
combination of $m_0,m_2,m_4,\dots,m_{2n}$.
\end{cor}


\section{Minkowski's question mark function}\label{sectmink?}

Given an irrational real number $x$ in $(0,1)$ with continued fraction 
expansion given by
$$x=[0;a_1,a_2,a_3,\dots]=\frac{1}{a_1+\frac{1}{a_2+\dots}},$$
Minkowski's question mark function is defined by
\begin{align}\label{def?}
?(x)&=-2\sum_{k=1}^\infty \frac{(-1)^k}{2^{a_1+\dots+a_k}}\ .
\end{align}

\begin{prop}\label{propConwMinkrec}
Minkowski's question mark function is an increasing
homeomorphism of $[0,1]$ such that $\CB\circ ?(x)=?\circ \CB(x)=x$.
\end{prop}
\proof[Proof
(given for the sake of self-containedness)] 
Since $\CB$ is an increasing 
homeomorphism of $[0,1]$,
it is enough to prove that $\CB\circ?(x)=x$ for every rational number $x$
in $[0,1]$.  We show this by induction on the length $n$ of
the continued fraction expansion $x=[0;a_1,a_2,\dots,a_n]$
of $x$.
The result clearly holds for $n=0$ (corresponding to $x=0$) and for $n=1$
(corresponding to the inverse of a non-zero natural integer).
Writing $\frac{p_k}{q_k}=[0;a_1,\dots,a_k]$ we have
$$?\left(\frac{p_{n-1}}{q_{n-1}}\right)=-2\sum_{j=1}^{n-1}\frac{(-1)^j}{2^{a_1+\dots+a_j}}=
\frac{m}{2^{a_1+\dots+a_{n-1}-1}}$$
for a suitable natural number $m$. We also have
$$?\left(\frac{p_{n-2}}{q_{n-2}}\right)=?\left(\frac{p_{n-1}}{q_{n-1}}\right)-
\frac{(-1)^n}{2^{a_1+
\dots+a_{n-1}-1}}.$$
Using  the induction hypothesis $\frac{p_k}{q_k}=\CB\circ ?
\left(\frac{p_k}{q_k}\right)$
for $k<n$ and applying (\ref{formCkab}) to 
$$?\left(\frac{p_n}{q_n}\right)=
\frac{m}{2^{a_1+\dots+a_{n-1}-1}}-\frac{(-1)^n}
{2^{a_n}2^{a_1+\dots+a_{n-1}-1}}$$
we get
$$\CB\circ?\left(\frac{p_n}{q_n}\right)
=\frac{a_np_{n-1}+p_{n-2}}{a_nq_{n-1}+q_{n-2}}=\frac{p_n}{q_n}\ .$$
\endproof
The graph of $?$ is well-known to behave in a self-similar way
as shown by the following well-known result:

\begin{prop} We have
\begin{align}\label{sym?}
?(1-x)&=1-?(x)
\end{align}
and 
\begin{align}\label{homogr?}
?\left(\frac{x}{1+x}\right)&=\frac{1}{2}?(x)
\end{align}
for all $x\in[0,1]$.
\end{prop}
\proof{} Identity (\ref{sym?}) follows from Proposition \ref{propsymC}
and Proposition \ref{propConwMinkrec}. Identity \ref{homogr?}
follows from the Definition (\ref{def?}) applied to 
$\frac{x}{1+x}=\frac{1}{1+1/x}=[0;1+a_1,a_2,a_3,\dots]$.
\endproof

The aim of this paper is to study the moments
$$m_z=\int_0^1 \CB(t)^zdt=\int_0^1x^zd?(x)$$
of the probability measure $d?$ associated to the distribution function
$?(x)=\int_0^xd?(t)$. The inequalities 
$$\sum_{n=1}^\infty \frac{1}{2^n}\left(1-\frac{1}{n}\right)^j\leq m_j\leq 
2\sum_{n=1}^\infty \frac{1}{2^n}\left(1-\frac{1}{n}\right)^j$$
coming from the evaluation $\CB(1-2^{-m})=1-\frac{1}{m+1}$, 
and the trivial upper bound $\left\vert {-z-j+1\choose j}\right\vert\leq 
{\vert z\vert+j-1\choose j}\leq (\vert z\vert+j)^{\vert z\vert}$
show that $z\mapsto m_z$ is an entire function of $\mathbb C$.

The function $m_z$ is also given by the expression
$$m_z=\sum_{k=0}^\infty {z+k-1\choose k}\gamma(k+z)m_k$$
(see (\ref{formIpa})) where $\gamma(z)$ is the entire function defined by
$$\gamma(z)=\sum_{n=1}^\infty \frac{1}{2^n}\frac{1}{(1+n)^z}\ .$$

We give the series expansion of the entire function $z\mapsto m_z$
at $z=0$ and study the asymptotics
of $m_z$ for real $z\longrightarrow \pm \infty$.

\begin{prop}\label{propidmzmmz}
We have the identities
$$m_z=\sum_{j=0}^\infty {z\choose j}(-1)^jm_j=\sum_{j=0}^\infty 
{-z+j-1\choose j}m_j$$
where ${z\choose j}=\frac{z(z-1)(z-2)\cdots (z-j+1)}{j!}$.
\end{prop}

The main contribution to $m_{-k}$ given by Proposition \ref{propidmzmmz}
corresponds to indices $j$ such that 
$\frac{k+j}{j}e^{-\sqrt{(\log 2)/j}}\sim 1$ yielding $j\sim\frac{k^2}{\log 2}$.

Thus we have for example
\begin{align*}
m_{-1}&=\sum_{n=0}^\infty m_n\\
m_{-2}&=\sum_{n=0}^\infty (n+1)m_n\\ 
m_{-3}&=\frac{1}{2}\sum_{n=0}^\infty (n+1)(n+2)m_n
\end{align*}
and more generally
$$m_{-n}=\sum_{k=0}^\infty {k+n-1\choose n-1}m_k\ .$$

\proof[Proof of Proposition \ref{propidmzmmz}] Proposition 
\ref{propsymC} implies the equalities
$$m_z=\int_0^1\CB(t)^zdt=\int_0^1(1-\CB(t))^zdt=
\sum_{j=0}^\infty {z\choose j}(-1)^jm_j=\sum_{j=0}^\infty 
{-z+j-1\choose j}m_j$$
which hold for all $z\in\mathbb C$ since $\CB(t)\in(0,1)$ for $t\in (0,1)$.
\endproof


\section{A few useful identities}\label{sectusefulid}

Almost all results of this paper are based on a few trivial 
identities, recorded in this Section for later use.

\subsection{Binomial coefficients}

\begin{lem} \label{lembinom}
We have the series expansion
\begin{align}\label{formpuissbin}
\frac{1}{(1-x)^n}&=\sum_{k=0}^\infty {n+k-1\choose k}x^k=
 \sum_{k=0}^\infty {n+k-1\choose n-1}x^k
\end{align}
for $x$ in the open complex unit-disc.
\end{lem}

\proof{}
Apply the equality 
${-n\choose k}=(-1)^k{n+k-1\choose k}$
(where ${x\choose k}=\frac{x(x-1)\cdots (x-k+1)}{k!}$) to
Newton's identity $(1+(-x))^{-n}=
\sum_{k=0}^\infty {-n \choose k}(-x)^k$ or use induction on $n$.
\endproof

\begin{rem} Lemma \ref{lembinom} has the following nice combinatorial proof:
$\frac{1}{(1-x)^n}$ is the generating series for colouring Easter eggs with
$n$ different colours (or, equivalently, for the number of monomials
in $n$ commuting variables). The $k$-th coefficient is thus given by
${k+n-1\choose n-1}$.
\end{rem}

\begin{lem} \label{lemsumidentity} We have
\begin{align*}
\sum_{l=k}^j{j\choose l}{l\choose k}x^l=
{j\choose k}x^k(x+1)^{j-k}
\end{align*}
In particular, for $x=\frac{-1}{2}$ we get
$$\sum_{l=k}^j{j\choose l}{l\choose k}\frac{1}{(-2)^l}=
{j\choose k}\frac{(-1)^k}{2^j}$$
\end{lem}

\proof{} Compare the coefficients 
${j\choose l}{l\choose k}$ and ${j\choose k}{j-k\choose l-k}$ 
of $x^l$ of  both sides.
\endproof

\subsection{Identities for the Stern sequence}

We recall that the Stern sequence $s:\mathbb N\longrightarrow \mathbb N$
is recursively defined by $s(0)=0,s(1)=1,s(2n)=s(n)$ and $s(2n+1)=s(n)+
s(n+1)$ for $n\geq 1$.

\begin{prop}\label{propsternids} For all $n\geq 0$ and for all $r$ 
such that $0\leq r\leq 2^n$,
the Stern sequence satisfies the identities
\begin{align}\label{sternidA}
s(2^n+r)&=s(2^n-r)+s(r),
\end{align}
\begin{align}\label{sternidB}
s(2^n+r)&=s(2^{n+1}-r),
\end{align}
\begin{align}\label{sternidC}
s(r)&=2s(2^n+r)-s(3\cdot 2^n+r).
\end{align}
\end{prop}

\proof{} 
The identities hold for $n=0$ and $r\in\{0,1\}$.
Since $s(2m)=s(m)$ they hold for $r$ even by induction.
For odd $r=2t+1< 2^{n+1}$, we sum the identities corresponding to 
$(n-1,t)$ and $(n-1,t+1)$ which hold by induction. 
The definition $s(2m+1)=s(m)+s(m+1)$ and induction 
implies the identities for 
odd $r$.
\endproof

The main idea of this paper is to apply Lemma 
\ref{lembinom} to the trivial identities
\begin{align}\label{idtrivsSa}
\frac{\alpha s+\beta S}{\gamma s+\delta S}
&=\frac{\alpha}{\delta}\frac{s}{S}\frac{1}{\left(1+\frac{\gamma}{\delta}
\frac{s}{S}\right)}+\frac{\beta}{\delta}\frac{1}{\left(1+\frac{\gamma}{\delta}\frac{s}
{S}\right)},
\end{align}
\begin{align}\label{idtrivsSb}
\frac{\alpha s+\beta S}{\gamma s+\delta S}
&=\frac{\alpha}{\gamma}+\left(\frac{\beta}{\delta}-\frac{\alpha}{\gamma}\right)
\frac{1}{\left(1+\frac{\gamma}{\delta}\frac{s}{S}\right)},
\end{align}
\begin{align}\label{idtrivsSc}
\frac{\alpha s+\beta S}{\gamma s+\delta S},
&=\frac{\beta}{\delta}+\left(\frac{\alpha}{\delta}-\frac{\beta\gamma}{\delta^2}\right)
\frac{s}{S}\frac{1}{\left(1+\frac{\gamma}{\delta}\frac{s}{S}\right)}.
\end{align}



\section{A simple set of linear equations for $m_{\mathbb N}$}
\label{sectsimpleset}

\begin{thm} \label{thmformuleA} The sequence
$m_{0}=1,m_{1}=\frac{1}{2},m_2,\dots$ of moments defined by $m_n=\int_0^1\CB(t)^ndt$ (see (\ref{defmoments}))
satisfies the equalities
\begin{align}\label{formthmA}
m_n
&=\frac{1}{2^n}\sum_{k=0}^{\lfloor n/2\rfloor}
{n\choose 2k}\phi_{2k}
\end{align}
where 
\begin{align}\label{formphi}
\phi_n=\sum_{k=0}^\infty {n+k-1\choose k}\frac{m_{n+k}}{2^{n+k}}\ .
\end{align}
\end{thm}

\begin{rem}\label{remprinccontr}
Since the increasing function
$$k\longmapsto \frac{{n+k-1\choose k}}{2^{n+k}}\frac{2^{n+k+1}}
{{n+k\choose k+1}}=2\frac{k+1}{n+k}$$
(for $k>0$ and $n$ a fixed natural integer) equals $1$ for $k=n-2$
and since the moments $m_n$ are slowly decreasing,
the main contribution to $\phi_n$ comes asymptotically from summands
with indices $k$ roughly equal to $n$.

The main contribution to $\phi_n$ is thus given by moments of the form
$m_{2n+l}$ with $l$ an element of $\mathbb Z$ of small absolute value.

Similarly, the main contribution to $m_n$ in Formula 
(\ref{formthmA}) corresponds asymptotically to indices $k\sim n/4$.
and involves thus mainly moments of the form $m_{n+l}$
for $l$ a small integer.
\end{rem}

Theorem \ref{thmformuleA} is an immediate consequence of the 
following result.

\begin{prop} \label{propforthmA}
For all $n\in\mathbb N$ we have the identities
$$\int_0^{1/2}\CB(t)^ndt=\frac{1}{2^{n+1}}\sum_{k=0}^n{n\choose k}(-1)^k\phi_k$$
and
$$\int_{1/2}^1\CB(t)^ndt=\frac{1}{2^{n+1}}\sum_{k=0}^n{n\choose k}\phi_k$$
with $\phi_k$ defined by Formula (\ref{formphi}).
\end{prop}

\begin{lem}\label{lemphi} We have
\begin{align}\label{formulaphin}
\phi_n&=2\int_{1/2}^1(2\CB(t)-1)^ndt\ .
\end{align}
for $\phi_n$ defined by Formula (\ref{formphi}).
\end{lem}

Corollary \ref{corsymmetric} shows that Lemma \ref{lemphi}
can be restated as 
$\phi_n=\int_0^1\left\vert 2\CB(t)-1\right\vert^ndt$.

\proof[Proof of Lemma \ref{lemphi}] 
Proposition \ref{propfond} and the definition of Riemann sums show that 
we have
\begin{align*}
2\int_{1/2}^1(2\CB(t)-1)^ndt&=2\lim_{l\rightarrow\infty}
\frac{1}{2^{l+1}}\sum_{r=0}^{2^l}
\left(2\frac{s(2^l+r)}{s(2^{l+1}+2^l+r)}-1\right)^n\\
&=\lim_{l\rightarrow\infty}
\frac{1}{2^l}\sum_{r=0}^{2^l}\left(\frac{2s(2^l+r)-s(3\cdot 2^l+r)}
{s(3\cdot 2^l+r)}\right)^n\ .
\end{align*}
Using  (\ref{sternidC}) we get
\begin{align}\label{calphinid}
2\int_{1/2}^1(2\CB(t)-1)^ndt&=\lim_{l\rightarrow\infty}
\frac{1}{2^l}\sum_{r=0}^{2^l}\left(\frac{s(r)}{2s(2^l+r)-s(r)}\right)^n
\end{align}
or equivalently
\begin{align*}
2\int_{1/2}^1(2\CB(t)-1)^ndt&=\lim_{l\rightarrow\infty}
\frac{1}{2^l}\sum_{r=0}^{2^l}\left(\frac{s(r)}{2s(2^l+r)}\right)^n
\left(\frac{1}{1-\frac{s(r)}{2s(2^l+r)}}\right)^n\ .
\end{align*}
Applying (\ref{formpuissbin}) we have
\begin{align*}
2\int_{1/2}^1(2\CB(t)-1)^ndt&=\sum_{j=0}^\infty{n+k-1\choose k}\frac{m_{n+k}}{2^{n+k}}
\end{align*}
which ends the proof.
\endproof

\proof[Proof of Proposition \ref{propforthmA}] 
Using Lemma \ref{lemphi}
we have
\begin{align*}
\frac{1}{2^{n+1}}\sum_{k=0}^n{n\choose k}(-1)^k\phi_k&=
\frac{1}{2^n}\int_{1/2}^1\sum_{k=0}^n{n\choose k}(1-2\CB(t))^kdt\\
&=\int_{1/2}^1(1-\CB(t))^ndt
\end{align*}
which equals 
$\int_0^{1/2}\CB(t)^ndt$ by (\ref{Csymmetric}).
This proves the first equality. 

The proof of the second equality is
similar and left to the reader.
\endproof


\subsection{Spectral properties}\label{subsectspect}

Theorem \ref{thmformuleA} expresses the 
moment-vector $(m_0,m_1,m_2,\dots)$ as a fixed point of a
continuous linear operator $T$ acting on the vector space $l^\infty(\mathbb R)$
of real bounded sequences. We study here 
a few spectral properties of $T$. They imply in particular uniqueness
of the fixed point
$(m_0,m_1,\dots)$ satisfying $m_0=1$.

We denote by $l^\infty=l^\infty(\mathbb R)$ the real Banach space
of bounded sequences with norm
$\parallel v\parallel_\infty=\sup_{n\in\mathbb N} \vert v_n\vert$ for 
$v=(v_0,v_1,\dots)$ in $l^\infty$. We set
$$\parallel U\parallel=\sup_{v\in l^\infty,\ \parallel v\parallel_\infty=1}
\parallel U(v)\parallel$$
for the norm $\parallel U\parallel$ of an endomorphism
$U\in\mathrm{End}(l^\infty)$. Similarly, we consider the norm
$$\parallel L\parallel=\sup_{v\in l^\infty,\ \parallel v\parallel_\infty=1}
\vert L(v)\vert$$
of a continuous linear form $L:l^\infty\longrightarrow \mathbb R$.

Formulae (\ref{formthmA}) and (\ref{formphi}) suggest to consider the 
sequence of operators
\begin{align}\label{formTi}
v&=(v_0,v_1,\dots)\longmapsto
T_n(v)=\frac{1}{2^n}\sum_{k=0}^{\lfloor n/2\rfloor}
{n\choose 2k}\sum_{l=0}^\infty {2k+l-1\choose l}\frac{v_{2k+l}}{2^{2k+l}}\ .
\end{align}

\begin{prop} \label{propnormT} Formula (\ref{formTi}) defines continuous
linear forms $T_0,T_1,T_2,\dots$ of norm $\parallel T_0\parallel =1$ and 
$\parallel T_n\parallel =\frac{1}{2}$ for $n\geq 1$.
\end{prop}

We define an endomorphism $T:l^\infty\longrightarrow l^{\infty}$ of the 
vector-space $l^\infty$ by setting $T=(T_0,T_1,T_2,\dots)$.
Proposition \ref{propnormT} and $T_0(v_0,v_1,\dots)=v_0$
imply the following result:

\begin{cor}\label{corattrfxpt} 
The restriction of the linear operator $T=(T_0,T_1,T_2,\dots)$ 
to the subspace
$l_0^\infty$ formed by all bounded sequences $(v_0,v_1,v_2,\dots)$
starting with $v_0=0$ yields an endomorphism of $l_0^\infty$ 
whose spectrum is contained in $\{z\in\mathbb C\ \vert
\ \vert z\vert\leq \frac{1}{2}\}$.

In particular, the linear map
$$v\longmapsto T(v)=(T_0(v),T_1(v),\dots)$$
defines a bounded linear operator of $l^\infty$ which
has a unique eigenvector of eigenvalue $1$
of the form  $\left(1,\frac{1}{2},\dots\right)$.
\end{cor}

The coordinates $(m_0,m_1,\dots)=\left(1,\frac{1}{2},\dots\right)$
of the unique eigenvector of eigenvalue $1$ of $T$ are of course the moments
$m_n=\int_0^1x^nd?x$ of the density function 
associated to Minkowski's question-mark function $?$.

\proof[Proof of Proposition \ref{propnormT}]
For $v\in l^\infty$ such that 
$\parallel v\parallel_\infty \leq 1$, we have
\begin{align*}
\vert T_n(v)\vert&\leq\frac{1}{2^n}\sum_{k=0}^{\lfloor n/2\rfloor}
{n\choose 2k}\frac{1}{2^{2k}}\sum_{j=0}^\infty{2k-1+j\choose j}\frac{1}{2^j}
\end{align*}
with equality if and only if $v$ is (up to a sign) the vector 
${\mathbf 1}=(1,1,1,\dots)$ with all coefficients equal to $1$.

Applying (\ref{formpuissbin}) we have thus
\begin{align*}
\parallel T_n\parallel&=\vert T_n({\mathbf 1})\vert\\
&=\frac{1}{2^n}\sum_{k=0}^{\lfloor n/2\rfloor}
{n\choose 2k}\frac{1}{2^{2k}}\left(\frac{1}{1-\frac{1}{2}}\right)^{2k}\\
&=\frac{1}{2^n}\sum_{k=0}^{\lfloor n/2\rfloor} {n\choose 2k}\\
&=\frac{1}{2^n}\left(\frac{(1+1)^n+(1-1)^n}{2}\right)\\
&=\left\lbrace\begin{array}{ll}
1&\hbox{if }n=0\\
\frac{1}{2}&\hbox{if }n\geq 1\end{array}\right.
\end{align*}
which completes the proof.
\endproof

\begin{rem}
Laplace's method shows that the coefficient
$$\frac{1}{2^n}\sum_{k=1}^{\lfloor n/2\rfloor}{n\choose 2k}
{m-1\choose 2k-1}\frac{1}{2^m}$$
of $v_m$ in $T_n$ given by Formula (\ref{formTi})
is asymptotically equal to 
$$\frac{1}{2}\frac{1}{\sqrt{2\pi n\mu(1+\mu)}}\left(
\frac{((1+\mu)/2)^{1+\mu}}{\mu^\mu}\right)^n$$
for $\mu=\frac{m}{n}$ having a bounded logarithm.
This coefficient is asymptotically maximal for $\mu=1$ and decays
exponentially fast otherwise. We have 
$$\lim_{n\rightarrow \infty}n\int_{0}^\infty 
\frac{1}{2}\frac{1}{\sqrt{2\pi n\mu(1+\mu)}}\left(
\frac{((1+\mu)/2)^{1+\mu}}{\mu^\mu}\right)^nd\mu=\frac{1}{2}$$
in agreement with Proposition \ref{propnormT}.
\end{rem}

\begin{rem} The linear operator $T$ has an unbounded eigenvector
of eigenvalue $\frac{1}{2}$ given by $w=(0,0,2,3,4,5,6,7,\dots)$
as can be seen as follows:
We have $T_0(w)=T_1(w)=0$. For $n\geq 2$, Formula (\ref{formTi}) with
$w=(0,0,2,3,4,5,6,\dots)$ boils down to 
$$T_n(w)=\frac{1}{2^n}\sum_{k=1}^{\lfloor n/2\rfloor}{n\choose 2k}
\sum_{l=0}^\infty{2k+l-1\choose l}\frac{2k+l}{2^{2k+l}}.$$

Computing the derivative $2k\frac{x^{2k-1}}{(1-x)^{2k+1}}$ of $\left(\frac{x}{1-x}\right)^{2k}$ at $x=\frac{1}{2}$ either directly or using the series expansion 
(\ref{formpuissbin}) given by Lemma \ref{lembinom} we get the identity 
\begin{align*}
4k&=\sum_{l=0}^\infty{2k+l-1\choose l}\frac{2k+l}{2^{2k+l}}\ .
\end{align*}
For $n\geq 2$ we have thus
\begin{align*}
T_n(w)&=\frac{1}{2^n}\sum_{k=0}^{\lfloor n/2\rfloor}{n\choose 2k}4k\\
&=\frac{1}{2^n}\left((1+x)^n+(1-x)^n\right)'\big\vert_{x=1}\\
&=\frac{n}{2}.\\
\end{align*}
\end{rem}

\subsection{Computational aspects}\label{sectcompaspects}

Theorem \ref {thmformuleA} is useful for computing numerical approximations 
of the first moments $m_0,m_1,\dots,m_N$ of Minkowski's question mark function. 

This can be done by computing an approximation $(\tilde m_0,\tilde m_1,
\dots,\tilde m_N)$ of the unique attracting 
fixed point $(\tilde m_0,\tilde m_1,
\dots)$ of the form $(1,\dots)$
of the linear operator $T\circ \pi_N$
where $\pi_N:l^\infty\longrightarrow l^\infty$ 
is the projection defined by
$$\pi_N(x_0,x_1,\dots,x_N,x_{N+1},\dots)=(x_0,x_1,\dots,x_N,0,0,0,\dots).$$

The error $\vert \tilde m_i-m_i\vert$ is of order $O(m_{N+1})=
O\left(N^{1/4}e^{-2\sqrt{N\log 2}}\right)$, see Formula (\ref{approxmnlambda}). 

Since the distance to the fixed point is essentially divided by 
$2$ under each iteration of 
$T\circ \pi_N$, the complexity of the resulting algorithm is
roughly of order $O\left(\sqrt{N/\log 2}N^2\right)$ if aiming
at maximal accuracy.

More precisely, the algorithm can be implemented as follows:

\noindent 010\quad $\tilde m_0:=1$,

\noindent 020\quad For $n=1,2,3,\dots,N$ do:

\noindent 030\quad\quad $\tilde m_n:=0$,

\noindent 040\quad End of loop over $n$,

\noindent 050\quad Iterate the following loop:

\noindent 060\quad\quad For $n=0,2,4,\dots,2\lfloor N/2\rfloor$ do:

\noindent 070\quad\quad\quad $b:=\frac{1}{2^n}$,

\noindent 080\quad\quad\quad $\tilde\phi_n:=0$, 

\noindent 090\quad\quad\quad For $k=0,1,2,\dots,N-n$
do: 

\noindent 100\quad\quad\quad\quad
$\tilde\phi_n:=\tilde\phi_n+b\tilde m_{n+k}$,

\noindent 110\quad\quad\quad\quad $b:=\frac{n+k}{2(k+1)}b$,

\noindent 120\quad\quad\quad End of loop over $k$,

\noindent 130\quad\quad End of loop over $n$,

\noindent 140\quad\quad For $n=1,2,3,\dots,N$ do:

\noindent 150\quad\quad\quad $b:=\frac{1}{2^n}$,

\noindent 160\quad\quad\quad $\tilde m_n:=0$,

\noindent 170\quad\quad\quad For $k=0,1,2,\dots,\lfloor N/2\rfloor$
do:

\noindent 180\quad\quad\quad\quad
$\tilde m_n:=\tilde m_n+b\tilde\phi_{2k}$,

\noindent 190\quad\quad\quad\quad $b:=\frac{(n-2k-1)(n-2k)}{(2k+1)(2k+2)}b$,

\noindent 200\quad\quad\quad End of loop over $k$,

\noindent 210\quad\quad End of loop over $n$,

\noindent 220\quad End of outer loop (starting at 050).

{\bf Comments:}\begin{enumerate}
\item{} Computations should be done over the real numbers 
with sufficient accuracy
(maximal achievable accuracy is of order $O(m_{N+1})$, see Section 
\ref{sectasympt} for estimations).

\item{} The range and increment of the loop-variable $n$ in line 060
is due to the fact that $m_1,\dots,m_N$ depend only on
$\phi_0,\phi_2,\phi_4,\dots,\phi_{2\lfloor N/2\rfloor}$ 
in Formula (\ref{formthmA}).

\item{} Instructions 070 and 150
need a loop in many programming languages.

\item{} The variable $b$ in line 070, 100, 110 corresponds to the factor
${n+k-1\choose k}\frac{1}{2^{n+k}}$ in Formula (\ref{formphi}).

\item{} The variable $b$ in line 150, 180,190 corresponds to the factor
$\frac{1}{2^n}{n\choose 2k}$ in Formula (\ref{formthmA}).

\item{} Maximal possible accuracy is achieved by iterating the 
outer loop (instructions 060-210) 
roughly $2\sqrt{N/\log 2}$ times, see Corollary \ref{corasympt}.

\item{} Using a known sequence of good 
approximations for $m_1,\dots,m_N$ instead of $0$ when initializing $\tilde m_1,\dots,\tilde m_N$ (instruction 030) decreases the number of useful (i.e. 
leading to significantly better precision) iterations for the outer loop.

\item{} A progressive increase of $N$ (starting from some small initial value) during 
the iteration of the outer loop yields a small speedup. 
\end{enumerate}

\section{Formulae of Alkauskas}\label{sectformAlk}

Theorem \ref{thmformuleA} is based on Riemann sums for the integral 
$$A=2\int_{1/2}^1(2\CB(t)-1)^ndt$$
obtained by subdividing 
the interval $\left[\frac{1}{2},1\right]$ into $2^l$ sub-intervals of 
equal length $\frac{1}{2^{l+1}}$.

In this section we give a new proof of some formulae obtained by Alkauskas 
by considering the infinite subdivision 
$$[0,1]=\{0\}\cup\dots \left[\frac{1}{2^h},\frac{1}{2^{h-1}}\right]\cup
\left[\frac{1}{2^{h-1}},\frac{1}{2^{h-2}}\right]\cup\dots\cup
\left[\frac{1}{4},\frac{1}{2}\right]\cup
\left[\frac{1}{2},\frac{1}{1}\right]$$
suggested by the easy evaluations
$\CB\left(\frac{1}{2^h}\right)=\frac{1}{h+1}$.

\begin{thm}\label{thmformAlk} We have 
\begin{align}\label{formIp}
m_n&=\sum_{h=1}^\infty \frac{1}{2^h}\frac{1}{(h+1)^n}
\sum_{k=0}^\infty {k+n-1\choose n-1}\frac{m_k}{(h+1)^k}
\end{align}
and 
\begin{align}\label{formInn}
m_n&=\sum_{h=1}^\infty \frac{1}{2^h}\frac{1}{h^n}
\sum_{k=0}^\infty {k+n-1\choose n-1}\frac{m_k}{(-h)^k}\ .
\end{align}
\end{thm}

\begin{rem} From a computational point of view it is perhaps useful 
to rewrite the formulae of Theorem \ref{thmformAlk} as
\begin{align}\label{formIpa}
m_n&=\sum_{k=0}^\infty {k+n-1\choose n-1}\gamma_{k+n}m_k
\end{align}
and 
\begin{align}\label{formIna}
m_n&=\sum_{k=0}^\infty {k+n-1\choose n-1}(-1)^kc_{k+n}
m_k
\end{align}
where 
$$\gamma_n=\sum_{k=1}^\infty \frac{1}{2^k(k+1)^n}=
2\mathrm{Li}_n\left(\frac{1}{2}\right)-1$$
and 
$$c_n=\sum_{k=1}^\infty \frac{1}{2^kk^n}=\mathrm{Li}_n\left(\frac{1}{2}\right)$$
where $\mathrm{Li}_n(x)=\sum_{k=1}^\infty \frac{x^k}{k^n}$
for $x$ in the open complex unit-disc.

Formula (\ref{formIp}) (or (\ref{formIpa})) should be preferred over 
(\ref{formInn}) (or (\ref{formIna})). It converges
faster (under iteration) and positivity of all coefficients ensures
numerical stability.

Precomputing (and storing) the constants $\gamma_k$ and using (\ref{formIpa})
needs only twice as much memory but provides a significant
speed-up.

Formula (\ref{formIna}) has been used by Alkauskas for numerical computations
of the first values of $m_n$, see Appendix 
A3 of \cite{AlMC} or Proposition 5 of \cite{AlGMJ}.

Since $\gamma_n\sim\frac{1}{2^{n+1}}$ for large $n$, the arguments of Remark
\ref{remprinccontr} show that the main contribution to $m_n$
in Formula (\ref{formIpa}) corresponds asymptotically to summands $k\sim n$
involving $m_{n-a},\dots,m_{n+a}$.
\end{rem}



\begin{prop}\label{propformIn} Setting
\begin{align}\label{defIhn}
I_h(n)&=\int_{2^{-h}}^{2^{-h+1}}\CB(t)^ndt\ .
\end{align} 
we have 
\begin{align}\label{formIn}
I_h(n)
&=\frac{1}{2^h}\sum_{k=0}^\infty {k+n-1\choose n-1}\frac{m_k}{(h+1)^{k+n}}
\end{align}
and 
\begin{align}
I_h(n)&=\frac{1}{2^h}\sum_{k=0}^\infty 
{k+n-1\choose n-1}(-1)^k\frac{m_k}{h^{k+n}}\ .
\end{align}
\end{prop}

\begin{lem}\label{lemcalChr} We have 
\begin{align}\label{idChrhl}
\CB\left(\frac{1}{2^h}+\frac{r}{2^{h+l}}\right)&=\frac{
s(2^l+r)}{(h+1)s(2^l+r)-s(r)}
\end{align}
and 
\begin{align}\label{idChrhla}
\CB\left(\frac{1}{2^h}+\frac{r}{2^{h+l}}\right)&=\frac{
s(2^l+r)}{hs(2^l+r)+s(2^l-r)}
\end{align}
for $0\leq r\leq 2^l$.
\end{lem}

\begin{rem} More generally, if 
$$\CB\left(\frac{q}{2^h}\right)=\frac{a}{b}\hbox{ and }
\CB\left(\frac{q+1}{2^h}\right)=\frac{c}{d}$$
with $(a,b)\in\mathbb N^2$ and  $(c,d)\in\mathbb N^2$
pairs of relatively prime natural numbers, then 
\begin{align*}
\CB\left(\frac{q}{2^h}+\frac{r}{2^{h+l}}\right)&=
\frac{as(2^l+r)+(c-a)s(r)}{bs(2^l+r)+(d-b)s(r)}\\
&=\frac{cs(2^l+r)+(a-c)s(2^l-r)}{ds(2^l+r)+(b-d)s(2^l-r)}
\end{align*}
for $l\in\mathbb N$ and for $r$ such that $0\leq r\leq 2^l$.
One can then apply (\ref{idtrivsSa}), (\ref{idtrivsSb}),
(\ref{idtrivsSc}) (or a similar identity) with $S=(2^l+r),\
s=s(r)$ in order to get Riemann sums for $\int_{a/b}^{c/d}
\CB(t)^ndt$.
\end{rem}

\proof[Proof of Lemma \ref{lemcalChr}] 
An induction on 
$h$ establishes the formula for $l=0$ (and $r\in\{0,1\}$).

An induction on $l$ (for constant $h$) ends the proof.
\endproof

\proof[Proof of Proposition \ref{propformIn}]
We have 
\begin{align*}
I_h(n)&=\frac{1}{2^h}
\lim_{l\rightarrow\infty}\frac{1}{2^l} \sum_{r=0}^{2^l}
\CB\left(\frac{1}{2^h}+\frac{r}{2^{h+l}}\right)^n\ .
\end{align*}
By (\ref{idChrhl}) we have
\begin{align*}
I_h(n)&=\frac{1}{2^h}
\lim_{l\rightarrow\infty}\frac{1}{2^l} \sum_{r=0}^{2^l}
\left(\frac{s(2^l+r)}{(h+1)s(2^l+r)-s(r)}\right)^n\\
&=\lim_{l\rightarrow\infty}\frac{1}{2^{l+h}}\sum_{r=0}^{2^l}
\frac{1}{(h+1)^n}
\left(\frac{1}{1-\frac{s(r)}{(h+1)s(2^l+r)}}\right)^n
\end{align*}
and (\ref{formpuissbin}) implies now
\begin{align*}
I_h(n)&=\lim_{l\rightarrow\infty}\frac{1}{2^{l+h}}\sum_{r=0}^{2^l}
\frac{1}{(h+1)^n}
\sum_{k=0}^\infty {k+n-1\choose n-1}\left(\frac{s(r)}{(h+1)s(2^l+r)}
\right)^k\\
&=\frac{1}{2^h}\sum_{k=0}^\infty {k+n-1\choose n-1}\frac{m_k}{(h+1)^{k+n}}\ .
\end{align*}
This proves the first equality.

The second equality follows from (\ref{sternidB}) applied
to (\ref{idChrhla}) yielding the identities
\begin{align*}
I_h(n)&=\frac{1}{2^h}\lim_{l\rightarrow\infty}\frac{1}{2^l} \sum_{r=0}^{2^l}
\left(\frac{s(2^l+r)}{hs(2^l+r)+s(r)}\right)^n\\
&=\lim_{l\rightarrow\infty}\frac{1}{2^{l+h}}\sum_{r=0}^{2^l}
\frac{1}{h^n}
\left(\frac{1}{1+\frac{s(r)}{hs(2^l+r)}}\right)^n\\
&=\frac{1}{2^h}\sum_{k=0}^\infty {k+n-1\choose n-1}(-1)^k
\frac{m_k}{h^{k+n}}\ .
\end{align*}
\endproof

\proof[Proof of Theorem \ref{thmformAlk}] 
Follows from 
$m_n=\sum_{h=1}^\infty I_h(n)$ where $I_h(n)$ is evaluated using Proposition
\ref{propformIn}.
\endproof


\section{Holomorphicity of $m_x$}\label{sectholom}

\begin{thm}\label{thmholom} (i) 
The map $n\longmapsto m_n$ extends to an entire function $x\longmapsto 
m_x$.

(ii) The series expansion of $x\longmapsto m_x$ at $x=0$ is given by
\begin{align}\label{formseriesmx}
\sum_{n=0}^\infty \frac{x^n}{n!}\sum_{k=n}^\infty c_{n,k}m_k
\end{align}
where
\begin{align}\label{formcnk}
\sum_{k=n}^\infty c_{n,k}x^k&=\left(\log(1-x)\right)^n\ .
\end{align}
Equivalently, the numbers $c_{n,k}$ are given by the equality
\begin{align}\label{cnkstirlingfirst}
c_{n,k}&=(-1)^k\frac{n!}{k!}s(k,n)
\end{align}
where the numbers $s(k,m)$ defined by
$\sum_{m=0}^ks(k,m)x^m=x(x-1)(x-2)\cdots(x-k+1)$ are Stirling numbers of
the first kind.
\end{thm}

\begin{rem} The rational numbers $c_{n,k}$ defined by (\ref{formcnk})
are given by the recursive formulae
$c_{0,0}=1,\ c_{0,k}=0$ if $k>0$ and 
$$c_{n+1,k}=-\sum_{j=1}^{k-1}\frac{c_{n,k-j}}{j}\ , n>0\ .$$
They are also defined by the equality
$$c_{n,k}=(-1)^n\sum_{a_1,\dots,a_n\geq 1,\ a_1+\dots+a_n=k}\frac{1}
{a_1\cdot a_2\cdots a_n}\ .$$
\end{rem}

\proof[Proof of Theorem \ref{thmholom}] 
Extending formula (\ref{defIhn})
by considering
$$I_h(x)=\int_{2^{-h}}^{2^{-h+1}}e^{x\log(\CB(t))}dt$$
for arbitrary $x\in\mathbb C$  (where $\log(\CB(t))\in \mathbb R$ denotes the 
usual logarithm of the strictly positive real number $\CB(t)$), the inequalities
$$\frac{1}{h+1}=\CB(2^{-h})\leq \CB(t)\leq \CB(2^{-h+1})=\frac{1}{h},\ t\in
[2^{-h},2^{-h+1}]$$ 
show
$$\vert I_h(x)\vert\leq \frac{1}{2^h}
\max_{t\in \left[\frac{1}{1+h},\frac{1}{h}\right]}\vert t^x\vert\leq
\frac{(1+h)^{\vert x\vert}}{2^h}\ .$$
This implies 
$$\left\vert
\sum_{h=1}^\infty I_h(x)\right
\vert\leq \sum_{h=1}^\infty  \frac{(1+h)^{\vert x\vert}}{2^h}
<\infty\ .$$
The map 
$x\longmapsto m_x=\sum_{h=1}^\infty I_h(x)$ defines thus an entire function which 
coincides with $m_x$ for $x\in \mathbb N$.

Using the symmetry $\CB(x)=1-\CB(x)$ we have
$$m_x=\lim_{l\rightarrow\infty}\frac{1}{2^l}\sum_{k=1}^{2^l}\left(\frac{s(k)}
{s(2^l+k)}\right)^x\\
=\lim_{l\rightarrow\infty}\frac{1}{2^l}\sum_{k=0}^{2^l-1}\left(1-\frac{s(k)}
{s(2^l+k)}\right)^x\ .$$
The $n-$th derivative of $m_x$ at $x=0$ evaluates thus to
\begin{align*}
&\lim_{l\rightarrow\infty} \frac{1}{2^l}\sum_{k=0}^{2^l-1}
\left(\log\left(1-\frac{s(k)}{s(2^l+k)}\right)\right)^n\\
&=\lim_{l\rightarrow\infty}\frac{1}{2^l}\sum_{k=0}^{2^l-1}
\left(-\sum_{j=1}^\infty\frac{1}{j}
\left(\frac{s(k)}{s(2^l+k)}\right)^j\right)^n
\end{align*}
which proves formula (\ref{formseriesmx}). 
\endproof

\begin{rem} Holomorphicity of $x\longmapsto m_x$ can also be proved using 
Proposition \ref{propidmzmmz}.
\end{rem}


\section{Two conjectural relations}\label{sectconjrel}

The derivative of the holomorphic function $x\longmapsto m_x$ (see
Theorem \ref{thmholom}) 
is given by 
$$-\sum_{n=1}^\infty \frac{m_n}{n}\sim 
-0.7924251285954891181912115152998913988894127820438$$
at the origin $x=0$. 
It coincides experimentally with the number 
$$-2\left(\log 2-\sum_{n=1}^\infty \frac{m_n}{n2^n}\right)$$
leading to the following conjectural identity.

\begin{conj}\label{conjderivmx} We have
\begin{align}\label{idconj1a}
\log 2&=\sum_{n=1}^\infty \frac{m_n}{2n}\left(1+\frac{1}{2^{n-1}}\right)\ .
\end{align}
\end{conj}

A variation is given by
\begin{conj}
\begin{align}\label{idconj1b}
\log 2&=\sum_{n=1}^\infty \frac{m_n}{n}\left(\frac{1}{2^n}-(-1)^n\right)\ .
\end{align}
\end{conj}


\section{A third set of formulae}\label{sectthirdset}

In this section we consider the partition
$$[0,1]\setminus\{1\}=\left[0,\frac{1}{2}\right]\cup\left[\frac{1}{2},\frac{3}{4}\right]\cup
\left[\frac{3}{4},\frac{7}{8}\right]\cup\left[\frac{7}{8},\frac{15}{16}\right]
\cup\dots.
$$
The resulting identities, well suited for computing asymptotics,
are given by the following result:

\begin{thm}\label{thmJh}
\begin{align*}
m_n&=\frac{1}{2}\sum_{j=0}^\infty{n+j-1\choose j}(-1)^jm_{n+j}\\
&+\sum_{h=2}^\infty \frac{1}{2^h}\left(\frac{h-1}{h}\right)^n
\sum_{k=0}^n{n\choose k}\left(\frac{1}{h(h-1)}\right)^k
\sum_{j=0}^\infty {k-1+j\choose j}(-1)^j\frac{m_{k+j}}{h^j}
\end{align*}
and
\begin{align*}
m_n&=
\sum_{h=1}^\infty\frac{1}{2^h}\left(\frac{h}{h+1}\right)^n
\sum_{k=0}^n{n\choose k}\left(\frac{-1}{h(h+1)}\right)^k
\sum_{j=0}^\infty {k-1+j\choose j}\frac{m_{k+j}}{(h+1)^j}\ .
\end{align*}
\end{thm}

\begin{rem} Only terms of order $h\sim\sqrt{n/\log 2}+O\left(n^{1/4}\right)$
yield large contributions to the first sum of the formulae in
Theorem \ref{thmJh}. Corresponding terms of the second sum 
(over $k$) for such contributions decay exponentially fast.
Terms of the third sum (over $j$) decay also exponentially fast for
fixed $h>1$ and for $k$ small.
\end{rem}

We set 
\begin{align}\label{defJn}
J_h(n)&=\int_{1-2^{-h+1}}^{1-2^{-h}}\CB(t)^ndt\ .
\end{align}

\begin{prop} \label{propJh} 
We have for all $h\in\mathbb N,\ h\geq 1$ the identities
\begin{align}\label{formJha}
J_h(n)&=\frac{1}{2^h}\left(\frac{h-1}{h}\right)^n
\sum_{k=0}^n {n\choose k}\left(\frac{1}{h(h-1)}\right)^k
\sum_{j=0}^\infty {k-1+j\choose j}(-1)^j\frac{m_{k+j}}{h^j}
\end{align}
and
\begin{align}\label{formJhb}
J_h(n)&=\frac{1}{2^h}\left(\frac{h}{h+1}\right)^n\sum_{k=0}^n {n\choose k}
\left(\frac{-1}{h(h+1)}\right)^k
\sum_{j=0}^\infty {k-1+j\choose j}\frac{m_{k+j}}{(h+1)^j}\ .
\end{align}
\end{prop}

Observe that (\ref{formJha}) boils down to 
\begin{align}\label{formJha0}
J_1(n)&=\frac{1}{2}\sum_{j=0}^\infty{n+j-1\choose j}(-1)^jm_{n+j}
\end{align}
for $h=1$.

\proof[Proof of Proposition \ref{propJh}]
Identity (\ref{sternidA}) of Proposition \ref{propsternids} implies
$$J_h(n)=\int_{2^{-h}}^{2^{-h+1}}\left(1-\CB(t)\right)^n dt\ .$$
Using Formula (\ref{idChrhla}) of Lemma \ref{lemcalChr} we get
\begin{align*}
J_h(n)&=\frac{1}{2^h}\lim_{l\rightarrow\infty}\frac{1}{2^l} \sum_{r=0}^{2^l}
\left(1-\frac{s(2^l+r)}{hs(2^l+r)+s(2^l-r)}\right)^n\\
&=\frac{1}{2^h}\lim_{l\rightarrow\infty}\frac{1}{2^l} \sum_{r=0}^{2^l}
\left(1-\frac{s(2^l+r)}{hs(2^l+r)+s(r)}\right)^n\\
&=\lim_{l\rightarrow\infty}\frac{1}{2^{h+l}} \sum_{r=0}^{2^l}
\left(\frac{h-1}{h}+\frac{1}{h^2}\left(\frac{s(r)}{s(2^l+r)+\frac{s(r)}{h}}\right)
\right)^n\\
&=\lim_{l\rightarrow\infty}\frac{1}{2^{h+l}}\left(\frac{h-1}{h}\right)^n\sum_{r=0}^{2^l}
\sum_{k=0}^n {n\choose k}\frac{1}{h^k(h-1)^k}
\left(\frac{s(r)}{s(2^l+r)+\frac{1}{h}s(r)}\right)^k\ .
\end{align*}
Using the identity
\begin{align*}
\left(\frac{s(r)}{s(2^l+r)+\frac{s(r)}{h}}\right)^k&=
\left(\frac{s(r)}{s(2^l+r)}\right)^k
\left(\frac{1}{1+\frac{1}{h}\frac{s(r)}{s(2^l+r)}}\right)^k\\
&=\sum_{j=0}^\infty
{k-1+j\choose j}\frac{(-1)^j}{h^j}\left(\frac{s(r)}{s(2^l+r)}
\right)^{k+j}
\end{align*}
obtained by applying formula (\ref{formpuissbin}), we get the first equation.

Starting with
\begin{align*}
J_h(n)&=\frac{1}{2^h}\lim_{l\rightarrow\infty}\frac{1}{2^l} \sum_{r=0}^{2^l}
\left(1-\frac{s(2^l+r)}{(h+1)s(2^l+r)-s(r)}\right)^n\\
&=\lim_{l\rightarrow\infty}\frac{1}{2^{h+l}} \sum_{r=0}^{2^l}
\left(\frac{h}{h+1}-\frac{1}{(h+1)^2}\left(
\frac{s(r)}{s(2^l+r)-\frac{s(r)}{h+1}}\right)\right)^n\\
&=\lim_{l\rightarrow\infty}\frac{1}{2^{h+l}}\left(\frac{h}{h+1}\right)^n\sum_{r=0}^{2^l}
\sum_{k=0}^n {n\choose k}
\frac{(-1)^k}{h^k(h+1)^k}\left(
\frac{s(r)}{s(2^l+r)-\frac{s(r)}{h+1}}\right)^k\ .
\end{align*}
and finishing as above yields the second identity.
\endproof

\proof[Proof of Theorem \ref{thmJh}]
Follows from Proposition \ref{propJh} applied
to the obvious identity $m_n=\sum_{h=1}^\infty J_h(n)$.
\endproof


\section{Asymptotics}\label{sectasympt}

We set 
\begin{align}\label{formulalambda}
\lambda&=\sum_{n=0}^\infty \frac{(\log 2)^n}{n!}m_n\ .
\end{align}
Numerically, $\lambda$ is approximately equal to
$$1.42815984554560290424313465212729430726822547802532544939052972\ .$$

\begin{thm}\label{thmasympt}
For every strictly positive $\epsilon$ there exists a natural integer $N$ 
such that 
$$\left|m_n-\lambda\sum_{h=2}^\infty \frac{1}{2^h}\left(\frac{h-1}{h}\right)^n\right|
\leq \epsilon\  m_n$$
if $n\geq N$.
\end{thm}

The error given by the asymptotic approximation
\begin{align}\label{asymptotsumapprmn}
m_n&\sim \lambda\sum_{h=2}^\infty \frac{1}{2^h}\left(1-\frac{1}{h}\right)^n
\end{align}
in Theorem \ref{thmasympt} is 
surprisingly small, see Section \ref{secterror}.

\begin{cor}\label{corasympt} 
We have 
\begin{align}\label{approxmnlambda}
m_n&\sim \lambda\frac{n^{1/4}}{(\log 2)^{3/4}}\sqrt{\frac{\pi}{2}}e^{-2\sqrt{n\log 2}}
\end{align}
for $n\rightarrow \infty$.
\end{cor}

Corollary \ref{corasympt} is of course equivalent to Theorem 1 in \cite{Alkasympt}. The constant $\lambda$ defined by (\ref{formulalambda})
is related to the constant 
$$c_0=\int_{0}^12^t\left(1-\frac{1}{2}?(t)\right)dt=\frac{1}{\log 2}-\frac{1}{2}
\int_0^1 2^t?(t)dt$$ 
in Theorem 1 of \cite{Alkasympt} by 
$$\lambda=c_02\log 2$$
and satisfies the following additional identities:
\begin{prop} \label{propvallambda} We have
$$\lambda=2\sum_{n=0}^\infty\frac{(-\log 2)^n}{n!}m_n
=\frac{4}{3}\sum_{n=0}^\infty \frac{(\log 2)^{2n}}{(2n)!}m_{2n}
=4\sum_{n=0}^\infty \frac{(\log 2)^{2n+1}}{(2n+1)!}m_{2n+1}\ .$$
\end{prop}

Observe that the constant $\lambda$ appears also in the asymptotic 
expression $\lambda\frac{n!}{(\log 2)^{n+1}}$ for $m_{-n}$, 
see \cite{AlGMJ} or Proposition \ref{propasympthmn}.

\begin{rem} A computation of $\lambda$ with high precision needs only 
relatively few initial values of $m_0,m_1,m_2,\dots$.
I ignore however a direct approach for accurately computing only the first few values of $m_2,m_3,m_4,\dots$.
\end{rem}

\begin{prop} \label{propmujasympt}
We have
$$\lim_{n\rightarrow\infty}
\left(\frac{n^{1/4}}{(\log 2)^{3/4}}\sqrt{\frac{\pi}{2}}e^{-2\sqrt{n\log 2}}
\right)^{-1}\left(\sum_{h=2}^\infty \frac{1}{2^h}\left(1-\frac{1}{h}\right)^n\right)
=1\ .$$
\end{prop}

\proof[Proof of Proposition \ref{propmujasympt}]
We apply Laplace's method to $\sum_{h=2}^\infty 
\frac{1}{2^h}\left(1-\frac{1}{h}\right)^n$.

The derivative
$$f_n'(x)=\frac{1}{2^x}\left(1-\frac{1}{x}\right)^n\frac{(n+x(1-x)\log 2)}{x(x-1)}$$
of the function $f_n(x)=\frac{1}{2^x}\left(1-\frac{1}{x}\right)^n$ 
has roots given by the solutions of $x^2-x=\frac{n}{\log 2}$.

Assuming $x$ real and positive, the positive root of $f_n'$ is 
given by 
$$\rho=\frac{1+\sqrt{1+4n/\log 2}}{2}= \sqrt{\frac{n}
{\log 2}}+\frac{1}{2}+\frac{1}{8}\sqrt{\frac{\log 2}{n}}+
O\left(\frac{1}{\sqrt{n}^3}\right)$$
and we have 
$$f_n(\rho)=\frac{1}{\sqrt{2}}e^{-2\sqrt{n\log 2}}
\left(1+O\left(\frac{1}{\sqrt{n}}\right)\right)\ .$$

A straightforward computation shows 
\begin{align}\label{formsecder}
f_n''(\rho)&=\frac{1}{\sqrt{2}}e^{-2\sqrt{n\log 2}}\left(-2\frac{
\sqrt{\log 2}^3}{\sqrt n}+O\left(\frac{1}{n}\right)\right)\ .
\end{align}
Applying Laplace's method
\begin{align*}
\int_2^\infty f_n(h)dh&\sim f_n(\rho)\int_{-\infty}^{\infty}e^{-\frac{f_n''(\rho)}
{f_n(\rho)}t^2/2}dt\\
&=\sqrt{\frac{2\pi f_n(\rho)^3}{-f_n''(\rho)}}
\end{align*}
to the integral approximation 
 $\int_2^\infty f_n(h)dh$ of $\sum_{h=2}^\infty f_n(h)$ we get the result.
\endproof

\begin{prop}\label{propexistenceA} For every $\epsilon>0$ 
there exists a natural integer $A$ such that 
$$0\leq m_n-\sum_{h\in\lfloor \sqrt{n/\log 2}-An^{1/4}\rfloor}^{\lfloor \sqrt{n/\log 2}
+An^{1/4}\rfloor}J_h(n)<\epsilon m_n$$
for all $n$ large enough with $J_h(n)=\int_{1-2^{-h+1}}^{1-2^{-h}}\CB(t)^ndt$
given by (\ref{defJn}).
\end{prop}

\proof
The easy evaluation
$\CB\left(1-\frac{1}{2^h}\right)=1-\frac{1}{h+1}$ for $h\in\mathbb N$ shows
$$1-\frac{1}{h}\leq \CB(t)\leq 1-\frac{1}{h+1}$$
for $t\in\left[1-\frac{1}{2^{h-1}},1-\frac{1}{2^h}\right]$ and we have
$$\frac{1}{2^h}\left(1-\frac{1}{h}\right)^x\leq 
J_h(x)\leq \frac{1}{2^h}\left(1-\frac{1}{h+1}\right)^x$$
for real positive $x$. Since the unique positive root of
the logarithmic derivative
$$\frac{df/dh}{f}=-\log 2+\frac{x}{h(h-1)}$$
with respect to $h$
of $f=\frac{1}{2^h}\left(1-\frac{1}{h}\right)^x$ is given by $h\sim\sqrt{x/\log 2}$
for large $x$, the decay of the function 
$$s\longrightarrow \frac{1}{2^{\sqrt{x/\log 2}+sx^{1/4}}}
\left(1-\frac{1}{\sqrt{x/\log 2}+sx^{1/4}}\right)^x$$ is exponentially fast
in $\vert s\vert$ for large $x$. This implies the result.
\endproof

\proof[Proof of Theorem \ref{thmasympt}] Setting 
$$\tilde J_h(n)=2^h\left(\frac{h}{h-1}\right)^nJ_h(n)$$
formula (\ref{formJha})
of Proposition \ref{propJh} shows the identities
\begin{align*}
\tilde J_h(n)&=\sum_{k=0}^n{n\choose k}
\frac{1}{h^k(h-1)^k}\sum_{j=0}^\infty {k-1+j\choose j}(-1)^j\frac{m_{k+j}}{h^j}\\
&=\sum_{k=0}^n\frac{(\log 2)^k}{k!}\frac{\prod_{j=1}^{k-1}\left(1-\frac{j}{n}\right)}
{\left(\frac{h(h-1)}{n}\log 2\right)^k}\sum_{j=0}^\infty {k-1+j\choose j}
(-1)^j\frac{m_{k+j}}{h^j}\ .
\end{align*}
For $k$ fixed and for $h=\sqrt{n/2\log 2}+O(n^{1/4})$
we have
$$\lim_{n\rightarrow\infty}\frac{\prod_{j=1}^{k-1}\left(1-\frac{j}{n}\right)}
{\left(\frac{h(h-1)}{n}\log 2\right)^k}=1$$
and we get the asymptotics
$$\tilde J_h(n)\sim
\sum_{k=0}^\infty\frac{(\log 2)^k}{k!}m_k=\lambda$$
for $h=\sqrt{n/\log 2}+O\left(n^{1/4}\right)$.

Proposition \ref{propexistenceA} shows now
\begin{align*}
m_n&\sim_{\epsilon}\sum_{h=\lfloor \sqrt{n/\log 2}-An^{1/4}\rfloor}^{
\lfloor \sqrt{n/\log 2}+An^{1/4}\rfloor}J_h(n)\\
&\sim_{\epsilon}\sum_{h=\lfloor \sqrt{n/\log 2}-An^{1/4}\rfloor}^{
\lfloor \sqrt{n/\log 2}+An^{1/4}\rfloor}
\frac{1}{2^h}\left(\frac{h-1}{h}\right)^n
\tilde J_h(n)\\
&\sim_{\epsilon}\lambda\sum_{h=2}^\infty\frac{1}{2^h}\left(\frac{h-1}{h}\right)^n
\end{align*}
for $n\rightarrow \infty$ and fixed $A$ (depending on $\epsilon$)
with $a\sim_{\epsilon} b$ denoting $\vert a-b\vert<\epsilon a$
for arbitrary small $\epsilon$ if $n$ is large enough.
\endproof

\proof[Proof of Proposition \ref{propvallambda}] Working with formula 
(\ref{formJhb}) we get the asymptotics
\begin{align*}
m_n&\sim\sum_{h=1}^\infty\frac{1}{2^h}\left(\frac{h}{h+1}\right)^n
\sum_{k=0}^\infty \frac{(-\log 2)^k}{k!}m_k\\
&=2\sum_{k=2}^\infty\frac{1}{2^h}\left(\frac{h-1}{h}\right)^n
\sum_{k=0}^\infty \frac{(-\log 2)^k}{k!}m_k
\end{align*}
which imply the first equality by comparing with Theorem \ref{thmasympt}.
The two other identities are easy consequences.
\endproof

\proof[Proof of Corollary \ref{corasympt}] Follows from Theorem \ref{thmasympt}
and Proposition \ref{propmujasympt}.
\endproof


\subsection{Asymptotic formula for $\phi_n$}
Using similar techniques, we get the asymptotic approximation
\begin{align}\label{asformphn}
\phi_n&\sim 2\lambda\sum_{h=3}^\infty \frac{1}{2^h}\left(1-\frac{2}{h}
\right)^n
\end{align}
(where $\lambda$ is given by (\ref{formulalambda})) for 
$\phi_n=2\int_{1/2}^1\left( 2\CB(t)-1\right)^ndt$,
see Formula (\ref{formulaphin}) in Lemma \ref{lemphi}.
The relative error seems again to be of order 
$O\left(\phi_n^{5/4}\right)$
and has again (suitably normalized) a more or less periodic
behaviour as a function of $\sqrt{n}$.

Using Laplace's method for the right side of (\ref{asformphn})
we get the simpler and less accurate expression
\begin{align}\label{asformphnsimpler}
\phi_n&\sim\lambda
\frac{(2n)^{1/4}\sqrt{\pi}}{\log(2)^{3/4}}e^{-2\sqrt{2n\log 2}}\ .
\end{align}


\subsection{A second asymptotic formula}\label{sectsecasympt}

The motivation for this section is the estimation of the order of the 
error in the asymptotic approximation (\ref{asymptotsumapprmn}).

A refinement of the Riemann sum underlying Formula 
(\ref{asymptotsumapprmn}) should yield
a slightly more accurate approximation for $m_n$.
The order of the difference between the two formulae 
should be a measure for the accuracy of (\ref{asymptotsumapprmn}). 

We subdivide the interval underlying
the integral $J_h(n)$ defined by (\ref{defJn}) into two intervals of equal
lengths. We have $J_h(n)=A_h(n)+B_h(n)$ where 
$$A_h(n)=\int_{1-2^{-h+1}}^{1-3\cdot 2^{-h-1}}\CB(t)^ndt\hbox{ and }
B_h(n)=\int_{1-3\cdot 2^{-h-1}}^{1-2^{-h}}\CB(t)^ndt\ .$$
We have 
\begin{align*}
A_h(n)&=\frac{1}{2^{h+1}}\lim_{l\rightarrow\infty}\frac{1}{2^l}
\sum_{r=0}^{2^l}\left(1-\frac{s(2^l+r)+s(r)}{hs(2^l+r)+(h+1)s(r)}\right)^n\\
&=\frac{1}{2^{h+1}}\lim_{l\rightarrow\infty}\frac{1}{2^l}
\sum_{r=0}^{2^l}\left(\frac{2h-1}{2h+1}-\frac{s(r)}{(2h+1)^2(s(2^l+r)-\frac{
h+1}{2h+1}s(r))}\right)^n\\
&=\frac{1}{2^{h+1}}\left(\frac{2h-1}{2h+1}\right)^n
\sum_{k=0}^n{n\choose k}\frac{(-1)^k}{(4h^2-1)^k}
\lim_{l\rightarrow\infty}\frac{2^k}{2^l}\sum_{r=0}^{2^l}\left(\frac{s(r)}{2s(2^l+r)-\frac{2h+2}{2h+1}s(r)}\right)^k
\end{align*}
which yields 
$$\lim_{h\rightarrow\infty}A_h(n)=\frac{1}{2^{h+1}}\left(\frac{2h-1}{2h+1}\right)^n
\sum_{k=0}^n{n\choose k}\frac{(-1)^k}{(4h^2-1)^k}2^k\phi_k$$
by Identity (\ref{calphinid}).

For $h=\sqrt{n/\log2}+O(n^{1/4})$ we have thus
$$A_h(n)\sim\frac{1}{2^{h+1}}\left(\frac{2h-1}{2h+1}\right)^n
\sum_{k=0}^\infty\frac{(-\log 2)^k}{2^k\ k!}\phi_k\ .$$
A similar calculation shows 
$$B_h(n)\sim\frac{1}{2^{h+1}}\left(\frac{2h-1}{2h+1}\right)^n
\sum_{k=0}^\infty\frac{(\log 2)^k}{2^k\ k!}\phi_k$$
for $h=\sqrt{n/\log2}+O(n^{1/4})$.

We get thus for large $n$ and $h=\sqrt{n/\log2}+O(n^{1/4})$ the approximation
$$J_h(n)\sim \frac{1}{2^h}\sum_{k=0}^\infty\frac{(\log 2)^{2k}}{2^{2k}\ (2k)!}
\phi_{2k}\ .$$

Setting
\begin{align}\label{definrho}
\rho&=\sum_{k=0}^\infty \frac{(\log 2)^{2k}}{2^{2k}\ (2k)!}
\phi_{2k}
\end{align}
we have asymptotically
$$
m_n\sim \rho\sum_{h=1}^\infty\frac{1}{2^h}\left(1-\frac{2}{2h+1}\right)^n\ .$$
Using Laplace's method we get the asymptotic approximation
$$\sum_{h=1}^\infty\frac{1}{2^h}\left(1-\frac{2}{2h+1}\right)^n\sim
\frac{n^{1/4}\sqrt{\pi}}{(\log 2)^{3/4}}e^{-2\sqrt{n\log 2}}.$$
This shows
\begin{align}\label{approxmnrho}
m_n&\sim\rho\frac{n^{1/4}\sqrt{\pi}}{(\log 2)^{3/4}}e^{-2\sqrt{n\log 2}}
\end{align}
and implies the identity
\begin{align}\label{idlambdarho}
\rho&=\frac{\lambda}{\sqrt{2}}
\end{align}
as can be seen by comparing the two asymptotic 
approximations (\ref{approxmnlambda}) and (\ref{approxmnrho}) of $m_n$.

The asymptotic formula 
\begin{align}\label{asymptformula1/2}
m_n&\sim\lambda\sum_{h=1}^\infty\frac{1}{2^{h+1/2}}\left(1-\frac{1}{h+1/2}\right)^n
\end{align}
should thus be slightly better than (\ref{asymptotsumapprmn}),
see Figure 1 in Section \ref{secterror}.


\subsection{An estimation for the error of the asymptotic formulae}\label{secterror}

Setting 
\begin{align}\label{defSxn}
x\longmapsto S_x(n)&=\sum_{h=1}^\infty \frac{1}{2^{h+x}}\left(1-\frac{1}{h+x}\right)^n,
\end{align}
the asymptotic formulae (\ref{asymptotsumapprmn}) and 
(\ref{asymptformula1/2}) can be rewritten as 
$m_n\sim \lambda S_0(n)$ and $m_n\sim \lambda S_{1/2}(n)$.
Since $x\longmapsto S_x(n)$ is almost $1$-periodic (for small positive
$x$ and huge fixed $n$) and oscillates experimentally 
around the exact value of the integral 
\begin{align}\label{defSint}
S_{\int}(n)&=\int_{1}^\infty\frac{1}{2^t}\left(1-\frac{1}{t}\right)^ndt,
\end{align}
it is tempting to rescale the errors $m_n-S_x(n)$ by the inverse of 
the factor
\begin{align}\label{formulakappan}
\kappa(n)&=\sqrt{\left(S_0(n)-S_{\int}(n)\right)^2+
\left(S_{1/4}(n)-S_{\int}(n)\right)^2}
\end{align}
given by the \lq\lq amplitude'' of the almost $1$-periodic function 
$x\longmapsto S_x(n)-S_{\int}(n)$.

The sequence $S_{\int}(n)$ of integrals
is easy to compute recursively: We have the initial values
\begin{align*}S_{\int}(0)&=
\int_1^\infty \frac{dt}{2^t}=\frac{1}{2\log 2}\\
S_{\int}(1)&=\int_1^\infty \frac{1}{2^t}\left(1-\frac{1}{t}\right)dt\\
&=\frac{1}{2\log 2}+\int_{-\infty}^{-\log 2}\frac{e^t}{t}dt?\\
&=\frac{1}{2\log 2}-\int_{\log 2}^\infty\frac{e^{-t}}{t}dt?\\
&=\frac{1}{2\log 2}+\mathrm{Ei}(-\log 2)\\
\end{align*}
(where $\mathrm{Ei}(x)=\int_{-\infty}^x\frac{e^t}{t}dt=-\int_{-x}^{\infty}\frac{e^{-t}}{t}dt$
is the exponential integral) and
integration by parts yields the recursion relation 
\begin{align}\label{recdefAn}
S_{\int}(n)&=\left(2+\frac{\log 2}{n-1}\right)S_{\int}(n-1)-S_{\int}(n-2).
\end{align}

The normalized errors 
\begin{align}
E_0(n)=&=\frac{1}{\kappa(n)}\left(m_n-\lambda\sum_{h=2}^\infty\frac{1}{2^{h}}\left(1-\frac{1}{h}\right)^n\right),\\
E_{1/2}(n)&=\frac{1}{\kappa(n)}\left(m_n-\lambda\sum_{h=1}^\infty\frac{1}{2^{h+1/2}}\left(1-\frac{1}{h+1/2}\right)^n\right),\\
E_{\int}(n)&=\frac{1}{\kappa(n)}\left(m_n-\lambda\int_{1}^{\infty}\frac{1}{2^t}\left(1-\frac{1}{t}\right)^ndt\right)
\end{align}
are depicted in Figure 1 representing the points $\left(\sqrt n,E_0(n)\right),
\left(\sqrt{n}, E_{1/2}(n)\right)$ and $\left(\sqrt{n},E_{\int}(n)\right)$
for $n$ in $\{100,\dots,400\}$. 
Points on the smallest sinusoidal curve are associated to $E_{\int}$,
points on the sinusoidal curve of intermediate size to $E_{1/2}$ and points
on the largest curve to $E_0$. In all three cases the error seems to be 
close to a damped periodic function of $\sqrt{n}$ of local amplitude
$O(\kappa(n))$. 
\begin{figure}[h]\label{asympterr}
\epsfxsize=8cm
\epsfysize=8cm
\epsfbox{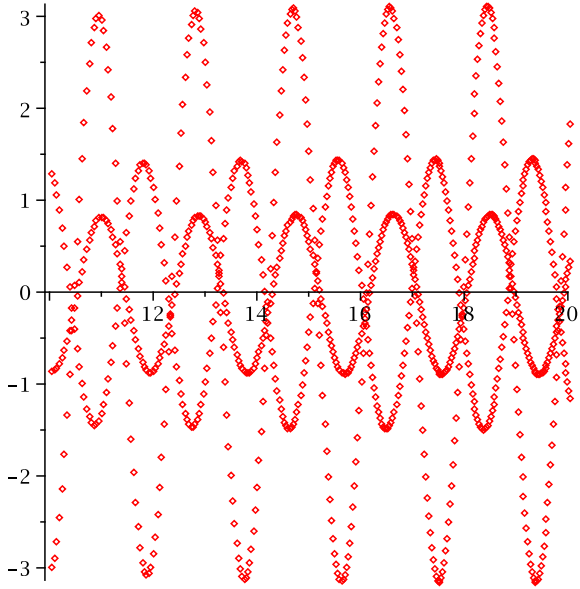}
\caption{Error of asymptotical formulae}
\end{figure}

\begin{rem} The existence of the linear recurrence relation (\ref{recdefAn})
implies the existence of asymptotic recurrence relations (given by the 
same formula) for the sequences $m_n$ and $S_x(n)$.

The asymptotic linear recurrence formula for $m_n$ can be improved 
into an affine asymptotic formula using ideas of the next Section.
\end{rem}

\begin{rem} It would be interesting to understand the asymptotic behaviour
of the amplitude $\kappa(n)$ given by Formula 
(\ref{formulakappan}).
(The number $\kappa(n)$ is essentially
the error term in Euler-MacLaurin's 
summation formula.) For moderate values of $n$ it seems to be comparable
to 
$$\frac{\sqrt{n}}{\log n\ \log\log n}e^{-9/2\sqrt{n\log 2}}$$
which implies $\kappa(n)<m_n^{9/4}$. The accuracy of 
the asymptotic formulae $m_n\sim \lambda S_*(n)$ (for $*=0,1/2$ and $\int$)
is thus surprisingly high.
\end{rem}

\subsection{Increasing accuracy}

The behaviour of the error-terms $E_*(n)$ occurring in the previous 
Section suggests to try an asymptotic formula of the form
\begin{align}\label{improvedasymptform}
m_n&\sim\lambda S_{\int}(n)+a\left(S_0(n)-S_{\int}(n)\right)
+b\left(S_{1/4}(n)-S_{\int}(n)\right)
\end{align}
with $\lambda$ defined by (\ref{formulalambda}) and 
$S_*(n)$ as in the previous Section.
Experimentally, such a formula seems to exist with 
\begin{align*}
a&\sim-.521901056340432536774725873446,\\
b&\sim-.148755851763595338634628933193.
\end{align*}
The term $\lambda S_{\int}(n)$ is of course the principal contribution
and plays the role of Formula (\ref{asymptotsumapprmn}) or (\ref{asymptformula1/2}).
The two remaining terms $a\left(S_0(n)-S_{\int}(n)\right)$ and 
$b\left(S_{1/4}(n)-S_{\int}(n)\right)$ sum up to a fairly 
regular (damped) oscillatory contribution of much lesser size.
More precisely, its local amplitude should be asymptotically equal to
$\frac{\lambda \sqrt{a^2+b^2}}{\kappa(n)}$ with
$\kappa(n)$ given by Formula (\ref{formulakappan}).

\subsection{An improved algorithm}

Accurate asymptotic approximations can be used for improving 
the algorithm given in Section \ref{sectcompaspects}.
Indeed, the cutoff at $N$ induces large relative errors for the last values of 
$\tilde \phi_n$.
It is thus natural to compute $\tilde \phi_0,\tilde \phi_2,\dots,
\tilde \phi_{2\lfloor N/2\rfloor}$ using the first $N+1$ values $\tilde m_0,\dots,\tilde m_N$
and $M-N$ additional values $\tilde m_{N+1},\dots,\tilde m_{M}$
given by asymptotic approximations of $m_{N+1},\dots,m_M$ 
(for $M>N$ a suitable integer
depending on $N$ and on the accuracy
of the chosen approximation).

We illustrate this by modifying the algorithm of Section 
\ref{sectcompaspects} using high-level instructions 
in order to involve the asymptotic 
approximation (\ref{asymptotsumapprmn}) (the approximation
(\ref{approxmnlambda}) is of much lesser interest):

Add the lines

\noindent 005\quad Precompute (and store) sufficiently accurate 
values $\tilde S(n)$ of $S_0(n)$ (or, slightly better, of $S_{\int}(n)$) for $n=N+1,\dots,M$.

\noindent 051\quad\quad Compute $\tilde \lambda:=\sum_{n=0}^M\frac{(\log 2)^n}{n!}\tilde m_n$,

\noindent 052\quad\quad Set $\tilde m_n=\lambda \tilde S(n)$ for 
$n=N+1,\dots,M$.

at the obvious locations.

Replace 090 by

\noindent 090\quad\quad\quad For $k=0,1,2,\dots,M-n$ do:

The resulting algorithm can easily be modified in order to work 
with other asymptotic approximations. The author used mainly
(\ref{improvedasymptform}) (this needs precomputations of
approximations for  $S_{\int}(n),S_0(n),S_{1/4}(n)$ with $n$ in
$\{N+1,\dots,M\}$).

Concerns using an algorithm based on a conjectural formula can 
be avoided by checking the final data
using a single iteration of (the main loop in) 
the original algorithm (described 
in Section \ref{sectcompaspects})
with a sufficiently high value $N'>N$ (with missing values
replaced by their (conjecturally very accurate) approximations).
The obtained data are exact up to an absolute error bounded
by $\max(\vert \epsilon\vert, m_{N'+1})$ with
 $\epsilon$ denoting the maximal modification of 
$\tilde m_1,\dots,\tilde m_{N'}$ during the final 
checking-run.

The improved version has smaller memory requirement and
a much better running time : 
The (conjectural) accuracy of the used approximation should more than
double the number of achievable correct digits for a given value of $N$.
In order to achieve the same accuracy, the original algorithm
has to be run with $N$ multiplied by more than $4$ which multiplies
the running time of the main loop by more than $16=4^2$.


\section{Values of $m$ at negative integers}\label{sectvalneg}

\begin{prop} \label{prophnhmn} 
The equality
\begin{align}\label{idnnegnpos}
m_{-n}&=m_n+\sum_{k=0}^{n-1}{n\choose k}\left(m_{-k}+m_{k}\right)
\end{align}
holds for $n\in\mathbb N$ a natural integer.
\end{prop}

\begin{rem} 
The generalization 
$$m_z=\frac{1}{2}\sum_{k=0}^\infty {-z\choose k}\left(m_{-k}+m_k\right)$$
of Proposition \ref{prophnhmn} fails 
for arbitrary complex values of $z$. Indeed, Proposition \ref{prophnhmn}
is based on the identity $(1+x)^z=\sum_{k=0}^\infty 
{z\choose k}x^k$ for arbitrary $x\in\mathbb R$ which breaks down if 
$-z$ is not in $\mathbb N$.
\end{rem}

\proof[Proof of Proposition \ref{prophnhmn}] We have for $n\in\mathbb N$
\begin{align*}
m_{-n}&=\lim_{l\rightarrow\infty}\frac{1}{2^l}\sum_{r=1}^{2^l}\left(
\frac{s(2^l+r)}{s(r)}\right)^n
\end{align*}
Using (\ref{sternidA}) we have
\begin{align*}
m_{-n}&=\lim_{l\rightarrow\infty}\frac{1}{2^l}\sum_{r=1}^{2^l}\left(
\frac{s(r)+s(2^l-r)}{s(r)}\right)^n\\
&=\lim_{l\rightarrow\infty}\frac{1}{2^l}\sum_{r=1}^{2^l}\left(1+
\frac{s(2^l-r)}{s(r)}\right)^n\\
&=\lim_{l\rightarrow\infty}\frac{1}{2^l}\left(\sum_{r=1}^{2^{l-1}}\left(1+
\frac{s(2^l-r)}{s(r)}\right)^n+\sum_{r=2^{l-1}+1}^{2^l}\left(1+
\frac{s(2^l-r)}{s(r)}\right)^n\right)\\
&=\lim_{l\rightarrow\infty}\frac{1}{2^l}\left(\sum_{r=1}^{2^{l-1}}\left(1+
\frac{s(2^l-r)}{s(r)}\right)^n+\sum_{r=0}^{2^{l-1}-1}\left(1+
\frac{s(r)}{s(2^l-r)}\right)^n\right)
\end{align*}
Using (\ref{sternidB}) we have thus
\begin{align*}
m_{-n}&=\lim_{l\rightarrow\infty}\frac{1}{2^l}\left(\sum_{r=1}^{2^{l-1}}
\left(1+\frac{s(2^{l-1}+r)}{s(r)}\right)^n+\sum_{r=1}^{2^{l-1}}\left(1+
\frac{s(r)}{s(2^{l-1}+r)}\right)^n\right)\\
&=\lim_{l\rightarrow\infty}\frac{1}{2^l}\sum_{k=0}^n{n\choose k}
\sum_{r=1}^{2^{l-1}}\left(
\left(\frac{s(2^{l-1}+r)}{s(r)}\right)^k+
\left(\frac{s(r)}{s(2^{l-1}+r)}\right)^k \right)\\
&=\frac{1}{2}\sum_{k=0}^n{n\choose k}\left(m_{-k}+m_k\right)
\end{align*}
which implies the result.
\endproof


\subsection{Matrices relating $m_{-\mathbb N}$ and $m_{\mathbb N}$}
\label{subsectmatrel-NN}

Identity (\ref{idnnegnpos}) of Proposition \ref{prophnhmn}
implies the existence of infinite 
lower diagonal triangular unipotent matrices 
$A,B=A^{-1}$ with integral coefficients such that 
$$\left(\begin{array}{c}m_0\\m_{-1}\\m_{-2}\\\vdots\end{array}\right)=
A\left(\begin{array}{c}m_0\\m_1\\m_2\\\vdots\end{array}\right)\hbox{ and }
\left(\begin{array}{c}m_0\\m_1\\m_2\\\vdots\end{array}\right)=
B\left(\begin{array}{c}m_0\\m_{-1}\\m_{-2}\\
\vdots\end{array}\right)\ .$$

The first few rows and columns of the matrices $A,B=A^{-1}$ are
$$\left(\begin{array}{rrrrrrrrrrr}
1\\
2&1\\6&4&1\\
26&18&6&1\\
150&104&36&8&1
\end{array}\right),\quad 
\left(\begin{array}{rrrrrrrrrrr}
1\\
-2&1\\
2&-4&1\\
-2&6&-6&1\\
2&-8&12&-8&1
\end{array}\right)
$$
and their coefficients are described by the following result.

\begin{prop} \label{propmatrixcoeffs} Let $\sigma_n,\ n\in\mathbb Z$ be a sequence (with values in a commutative ring containing $1$) indexed by
the set $\mathbb Z$ of all integers such that 
$$\sigma_{-n}-\sigma_n=\sum_{k=0}^{n-1}{n\choose k}
\left(\sigma_{-k}+\sigma_k\right)\ .$$

Then
\begin{align*}
\sigma_{-i}&=\sum_{j=0}^i \alpha_{i,j}\sigma_j\\
\sigma_i&=\sum_{j=0}^i \beta_{i,j}\sigma_{-j}
\end{align*}
for all $i\in\mathbb N$
where $\alpha_{i,j},\beta_{i,j},\ 0\leq i,j$ are integers given by the
formulae
\begin{align}\label{defalphij}
\alpha_{i,j}&={i\choose j}\sum_{h=1}^\infty \frac{h^{i-j}}{2^h}
\end{align}
and 
$$\beta_{i,j}=\left\lbrace\begin{array}{cl}
1&\hbox{if }i=j\\
2(-1)^{i+j}{i\choose j}\ &\hbox{otherwise.}\end{array}\right.$$
In particular, the matrices $A$ and $B$ with 
coefficients $\alpha_{i,j},\beta_{i,j},
\ 0\leq i,j$ are mutually inverse lower triangular unipotent 
integral matrices.
\end{prop}

\proof{}
We have $\alpha_{i,i}=1$ as required and the matrix $A$
is clearly lower triangular. The proof is now by 
induction on the row-index $i$ of the coefficients $\alpha_{i,j}$ for $A$.
Equation (\ref{idnnegnpos}) of Proposition \ref{prophnhmn} 
shows that we have
$$\alpha_{i+1,j}=\sum_{l=0}^i{i+1\choose l}\left(\alpha_{l,j}+\delta_{l,j}\right)$$
for $i+1\geq j$, 
where $\delta_{l,j}=1$ if $l=j$ and $\delta_{l,j}=0$ otherwise.

We get 
\begin{align*}
\alpha_{i+1,j}&={i+1\choose j} -{i+1\choose i+1}{i+1\choose j}\sum_{h=1}^\infty 
\frac{h^{i+1-j}}{2^h}\\
&\quad +\sum_{h=1}^\infty \sum_{k=0}^{i+1}
{i+1\choose k}{k\choose j}\frac{h^{k-j}}{2^h}\\
&={i+1\choose j}-{i+1\choose j}\sum_{h=1}^\infty \frac{h^{i+1-j}}{2^h}
+\sum_{h=1}^\infty \frac{{i+1\choose j}}{2^h}\sum_{k=0}^{i+1}
{i+1-j\choose k-j}h^{k-j}\\
&={i+1\choose j}\left(1-\sum_{h=1}^\infty \frac{h^{i+1-j}}{2^h}
+\sum_{h=1}^\infty \frac{(h+1)^{i+1-j}}{2^h}\right)\\
&={i+1\choose j}\sum_{h=1}^\infty \frac{h^{i+1-j}}{2^h}\ .
\end{align*}
This implies the formula for the coefficients of $A$.

We prove the formula for the coefficients of the inverse
matrix $B=A^{-1}$ by computing the product $AB$. We have
\begin{align*}
\sum_{k=j}^i\alpha_{i,k}\beta_{k,j}&=
2\sum_{k=j}^i\alpha_{i,k}(-1)^{k+j}{k\choose j}-\alpha_{i,j}\\
&=2\sum_{h=1}^\infty\sum_{k=j}^i\frac{h^{i-k}}{2^h}{i\choose k}{k\choose j}
(-1)^{k+j}-\alpha_{i,j}\\
&=\sum_{h=1}^\infty \frac{(-1)^jh^i}{2^{h-1}}\sum_{k=j}^i{i\choose k}{k\choose j}\frac{1}{(-h)^k}-\alpha_{i,j}\ .
\end{align*}
Identity \ref{lemsumidentity} of Lemma \ref{lemsumidentity} shows that 
this simplifies to
$$\sum_{h=1}^\infty\frac{(h-1)^{i-j}}{2^{h-1}}{i\choose j}
-\alpha_{i,j}\ .$$
This equals $1$ if $i=j$ and $0$ for $i>j$ by 
definition of $\alpha_{i,j}$.
\endproof

The sum appearing in (\ref{defalphij}) defines natural integers
having a recursive definition:
\begin{prop} \label{propcoeffsA}
The natural integers
$$\gamma_n=\sum_{h=1}^\infty \frac{h^n}{2^h}$$
(appearing in (\ref{defalphij})) 
have the recursive definition $\gamma_0=1$ and
$$\gamma_n=1+\sum_{j=0}^{n-1}{n\choose j}\gamma_j$$
for $n\geq 1$. 
\end{prop}

The sequence of integers $\gamma_0,\gamma_1,\dots$ starts as 
$$1,2,6,26,150,1082,9366,94586,1091670,\dots\ ,$$
see sequence A629 of \cite{OEIS}.

\proof[Proof of Proposition \ref{propcoeffsA}] 
We have
\begin{align*}
\gamma_n&=\frac{1}{2}+\frac{1}{2}\sum_{h=1}^\infty
\frac{(h+1)^n}{2^h}\\
&=\frac{1}{2}+\frac{1}{2}\sum_{j=0}^n{n\choose j}\sum_{h=1}^\infty
\frac{h^n}{2^h}\\
&=\frac{1}{2}+\frac{1}{2}\sum_{j=0}^n{n\choose j}\gamma_j
\end{align*}
which implies the result.
\endproof

\begin{rem} Lower triangular matrices with lower triangular coefficients
$\gamma_{i,j}$ of the form ${i\choose j}c_{i-j}$ for some sequence
$c_0,c_1,\dots$ form a commutative algebra. Indeed, the map
associating to such a matrix with coefficients ${i\choose j}c_{i-j}$
the formal exponential power series $\sum_{n=0}^\infty c_n\frac{x^n}{n!}$
defines an isomorphism of algebras onto the algebra of formal exponential
power series (with product given by the obvious \lq\lq bilinear\rq\rq
extension of $\frac{x^i}{i!}\frac{x^j}{j!}=
{i+j\choose i}\frac{x^{i+j}}{(i+j)!}$). 
The easy equality $\sum_{n=0}^\infty \beta_{n,0}\frac{x^n}{n!}=1+2\sum_{n=1}^\infty 
(-1)^n\frac{x^n}{n!}=2e^{-x}-1$
shows thus the identity
$$\sum_{n=0}^\infty \gamma_n\frac{x^n}{n!}=\sum_{n=0}^\infty a_{n,0}\frac{x^n}{n!}=\frac{1}{2e^{-x}-1}\ .$$
\end{rem}

Proposition \ref{prophnhmn} and Proposition \ref{propmatrixcoeffs} imply 
the following result:

\begin{cor}\label{cormmnsummk}
We have
\begin{align}\label{idmnneg}
m_{-n}=\sum_{k=0}^n {n\choose k}\gamma_{n-k}m_k
\end{align}
(with $\gamma_n$ defined by Proposition \ref{propcoeffsA}) and
\begin{align}\label{idmnpos}
m_n&=m_{-n}+2\sum_{k=0}^{n-1}(-1)^{n+k}{n\choose k}m_{-k}
\end{align}
for all $n$ in $\mathbb N$.
\end{cor}

Corollary \ref{cormmnsummk} is better suited than Proposition 
\ref{propidmzmmz} for computing values $m_{-\mathbb N}$
using $m_{\mathbb N}$. It involves only finitely many terms of $m_{\mathbb N}$
with coefficients which are decreasing.
(The main contribution to $m_{-n}$
given by the formula of Proposition \ref{propidmzmmz}
corresponds to summands indexed by integers close to $\frac{n^2}{\log 2}$.)

Combining Formula (\ref{idmnneg}) of Corollary \ref{cormmnsummk} with
Proposition \ref{propidmzmmz} we get:

\begin{cor}\label{corseqids} 
We have for all $n$ in $\mathbb N$ the identity
\begin{align}\label{seqidn}
\sum_{k=0}^n {n\choose k}\gamma_{n-k}m_k&
=\sum_{j=0}^\infty {n+j-1\choose j}m_j\ .
\end{align}
\end{cor}

Corollary \ref{corseqids} yields
\begin{align*}
2m_0+m_1&=\sum_{j=0}^\infty m_j,\\
6m_0+4m_1+m_2&=\sum_{j=0}^\infty (j+1)m_j,\\
26m_0+18m_1+6m_2+m_3&=\sum_{j=0}^\infty{j+2\choose 2}m_j,\\
150m_0+104m_1+36m_2+8m_3+m_4&=\sum_{j=0}^\infty{j+3\choose 3}m_j.
\end{align*}
Using the easy evaluation $m_0=1,m_1=\frac{1}{2}$, the case $n=1$ yields
the nice evaluation
\begin{align}\label{id5/2}
\sum_{j=0}^\infty m_j&=\frac{5}{2}
\end{align}
which can be used as an accuracy-check for numerical computations.

Similarly, using $n=2$, we get the identity $m_2+8=\sum_{j=0}^\infty (j+1)m_n$.
Subtraction of (\ref{id5/2}) yields
\begin{align}\label{identityn=2}
\sum_{j=1}^\infty j\ m_j&=m_2+\frac{11}{2}.
\end{align}


\subsection{Asymptotics for $m_{-n}$}

\begin{prop} \label{propasympthmn} We have 
$$\lim_{n\rightarrow\infty} 
m_{-n}\frac{(\log 2)^{n-1}}{n!}=\lambda$$
for
$\lambda=\sum_{k=0}^\infty \frac{(\log 2)^k}{k!}m_k$
given by (\ref{formulalambda}).
\end{prop}

The following easy result is probably well-known:

\begin{lem} \label{lemasymptkn2k} We have
$$\lim_{n\rightarrow\infty} \frac{(\log 2)^{n+1}}{n!}
\sum_{k=1}^\infty \frac{k^n}{2^k}=1\ .$$
\end{lem}

\proof[Proof of Lemma \ref{lemasymptkn2k}] 
We apply Laplace's method
to $\int_1^\infty x^n2^{-x}dx\sim\sum_{k=1}^\infty k^n2^{-k}$.

The derivative 
$$\frac{d}{dx}\left(x^ne^{-x\log 2}\right)=
\left(\frac{n}{x}-\log 2\right)x^ne^{-x\log 2}$$
of $x^n2^{-x}$
has a unique strictly positive root at $\frac{n}{\log 2}$ 
and second derivative
$-\frac{(\log 2)^2}{n}\left(\frac{n}{\log 2}\right)^ne^{-n}$ 
at the critical point $x=\frac{n}{\log 2}$ corresponding to the maximum
$\left(\frac{n}{\log 2}\right)^ne^{-n}$ of the function 
$x\longmapsto x^n2^{-x}$.

Laplace's method yields thus the asymptotics
\begin{align*}
\sum_{k=1}^\infty\frac{k^n}{2^k}&\sim
\left(\frac{n}{\log 2}
\right)^ne^{-n}\int_{-\infty}^\infty e^{-\frac{(\log 2)^2}{2n}x^2}dx\\
&=\sqrt{\frac{2\pi n}{(\log 2)^2}}\left(\frac{n}{\log 2}
\right)^ne^{-n}\\
&=\frac{1}{(\log 2)^{n+1}}\sqrt{2\pi n}\frac{n^n}{e^n}\\
&\sim \frac{n!}
{(\log 2)^{n+1}}
\end{align*}
where the last asymptotic equivalence follows from 
Stirling's formula $n!\sim\sqrt{2\pi n}\frac{n^n}{e^n}$.
\endproof

\proof[Proof of Proposition \ref{propasympthmn}]
Using Corollary \ref{cormmnsummk} 
and the asymptotics $\sum_{k=1}^\infty \frac{k^n}{2^k}\sim\frac{n!}
{(\log 2)^{n+1}}$ given by Lemma \ref{lemasymptkn2k}, 
we get the asymptotics
\begin{align*}
m_{-n}&=\sum_{k=0}^n\alpha_{n,k}m_k\sim \sum_{k=0}^n
{n\choose k}\frac{(n-k)!}{(\log 2)^{n-k+1}}m_k\\
&=\frac{n!}{(\log 2)^{n+1}}\sum_{k=0}^n\frac{(\log 2)^k}{k!}m_k\sim \lambda
\frac{n!}{(\log 2)^{n+1}}
\end{align*}
with $\lambda$ given by (\ref{formulalambda}).
\endproof

\section{Geometric means for the Stern sequence}\label{sectgeommeans}

It is an easy exercise to compute the arithmetic mean
$\frac{1}{2^n}\sum_{j=2^n}^{2^{n+1}}s(j)$.

The following result gives asymptotics for the geometric mean:

\begin{thm}\label{thmgeommean} There exists a real constant $\beta$ such that
$$\lim_{n\rightarrow\infty}e^{-(n\alpha+\beta) 2^n}\prod_{j=2^n}
^{2^{n+1}}s(j)=1$$
where
\begin{align*}\alpha&=
\log 2-\sum_{j=1}^\infty \frac{m_j}{j2^j}\\
&\sim 
.39621256429774455909560575764994569944470639102190\ .
\end{align*}
\end{thm}

\begin{rem} The constant $\alpha$ is involved in the Hausdorff dimension
of growth points for $?(x)$, see Kinney or Alkauskas. See also 
Conjecture \ref{conjderivmx}
for a conjectural manifestation of $\alpha$.

I am not aware of the existence of an efficient method 
for computing the value of $\beta\sim -.0851895$ with high precision.
\end{rem}

\begin{lem}\label{lemtrapezoidal} Given an increasing function $\varphi:[0,1]\longrightarrow
[0,1]$ and a strictly positive natural integer $N$ we have
$$\left\vert\int_0^1\varphi(t)dt-\frac{1}{N}\sum_{k=0}^{N-1}
\frac{\varphi\left(\frac{k}{N}\right)+\varphi\left(\frac{k+1}{N}\right)}{2}
\right\vert\leq\frac{1}{2N}\ .$$
\end{lem}

\proof{} The error of the trapezoidal rule 
$$\int_a^b\varphi(t)dt\sim(b-a)\frac{\varphi(a)+\varphi(b)}{2}$$
is bounded by $\left\vert\frac{(b-a)(\varphi(b)-\varphi(a))}{2}\right\vert$
if $\varphi$ is monotonous. 
\endproof

\proof[Proof of Theorem \ref{thmgeommean}]
We consider
\begin{align*}
S(n)&=\frac{1}{2^n}\sum_{k=2^n}^{2^{n+1}}\log(s(k))\\
&=\frac{1}{2^n}\sum_{k=1}^{2^n}\frac{\log(s(2^n+k-1))+\log(s(2^n+k))}{2}
\end{align*}
where the second identity follows from the evaluations $s(2^n)=s(2^{n+1})=1$.
Using (\ref{sternidB}) we get
\begin{align*}
S(n+1)&=\frac{1}{2^{n+1}}
\sum_{k=1}^{2^{n+1}}\frac{\log(s(2^{n+1}+k-1))+\log(s(2^{n+1}+k))}{2}\\
&=\frac{1}{2^{n+1}}
\sum_{k=1}^{2^{n}}\frac{\log(s(2^{n+1}+k-1))+\log(s(2^{n+1}+k))}{2}\\
&\quad +\frac{1}{2^{n+1}}
\sum_{k=1+2^n}^{2^{n+1}}\frac{\log(s(2^{n+1}+k-1))+\log(s(2^{n+1}+k))}{2}\\
&=\frac{2}{2^{n+1}}
\sum_{k=1}^{2^n}\frac{\log(s(3\cdot 2^n+k-1))+\log(s(3\cdot 2^n+k))}{2}
\end{align*}
and (\ref{sternidC}) yields
\begin{align*}
S(n+1)&=\frac{1}{2^n}
\sum_{k=1}^{2^n}\frac{\log(2s(2^n+k-1)-s(k-1))+\log(2s(2^n+k)-s(k))}{2}\\
&=\log(2)+\frac{1}{2^n}\sum_{k=1}^{2^n}\frac{
\log(s(2^n+k-1)-\frac{1}{2}s(k-1))+\log(s(2^n+k)-\frac{1}{2}s(k))}{2}
\end{align*}
Using $\frac{d^k}{dx^k}\log(u-vx)=-(k-1)!\frac{v^k}{(u-vx)^k}$ for $k\geq 1$,
we get
\begin{align*}
S(n+1)&=\log(2)+S(n)-\sum_{j=1}^\infty\frac{1}{j2^j}
\frac{1}{2^n}\sum_{k=1}^{2^n}\frac{1}{2}\left(
\left(\frac{s(k-1)}{s(2^n+k-1)}\right)^j+
\left(\frac{s(k)}{s(2^n+k)}\right)^j
\right)
\end{align*}
which implies 
$$\lim_{n\rightarrow\infty}(S(n+1)-S(n))=\log 2-\sum_{j=1}^\infty \frac{m_j}{j2^j}=\alpha$$
by Proposition \ref{propfond}.

Lemma \ref{lemtrapezoidal} shows
\begin{align*}
&\left\vert\sum_{n=0}^\infty\sum_{j=1}^\infty
\frac{1}{j2^j}\left(m_j-\frac{1}{2^n}\sum_{k=1}^{2^n}\frac{1}{2}\left(
\left(\frac{s(k-1)}{s(2^n+k-1)}\right)^j+
\left(\frac{s(k)}{s(2^n+k)}\right)^j
\right)\right)\right\vert\\
&\leq\sum_{n=0}^\infty\sum_{j=1}^\infty \frac{1}{2^j}\frac{1}{2^n}\leq
\sum_{n=0}^\infty\sum_{j=1}^\infty\frac{1}{2^{n+j}}=2\ .
\end{align*}
This proves the existence of $\beta$.
\endproof

\noindent Roland BACHER, 

\noindent Univ. Grenoble Alpes, Institut Fourier, 

\noindent F-38000 Grenoble, France.
\vskip0.5cm
\noindent e-mail: Roland.Bacher@ujf-grenoble.fr

\end{document}